\documentclass[letterpaper,11pt]{article} 

\usepackage{amsmath,amssymb} 
\usepackage{cite} 
\usepackage[margin=1in]{geometry} 
\usepackage{graphicx} 
\usepackage{hyperref} 
\usepackage{multirow} 
\usepackage{subfig} 

\begin{document}

\title{\textbf{Optimal Multi-Mode Propulsion Mission \\ Design Using Direct Collocation}}
\author{
George V. Haman III\thanks{Ph.D. Candidate, Department of Mechanical and Aerospace Engineering. Email: georgehaman@ufl.edu.}
~and Anil V. Rao\thanks{Professor, Department of Mechanical and Aerospace Engineering. Email: anilvrao@ufl.edu.} \vspace{12pt} \\ {\em{University of Florida}} \\ {\em{Gainesville, FL 32611}}
}
\date{}
\maketitle

\newcommand{\genericmathaccent}[2]{\vbox{\offinterlineskip\ialign{##\cr\hidewidth$\scriptstyle#1$\hidewidth\cr\noalign{\kern-.5ex}$#2$\cr}}} 
\newcommand{\genericmathaccenti}[2]{\genericmathaccent{\mkern1mu#1{}}{#2}} 
\newcommand{\genericmathaccentii}[2]{\genericmathaccent{\mkern1mu#1{}\mkern3mu#1{}}{#2}} 
\newcommand{\atantwo}{\text{atan2}} 
\newcommand{\dnu}{\ensuremath{\text{d}\nu}} 
\newcommand{\drv}[2]{\ensuremath{\dfrac{\text{d}{#1}}{\text{d}{#2}}}} 
\newcommand{\grad}[2]{\ensuremath{\nabla_{\mkern-3mu#1}{#2}}} 
\newcommand{\ildrv}[2]{\ensuremath{\mrm{d}{#1}/\mrm{d}{#2}}} 
\newcommand{\ilndrv}[3]{\ensuremath{\mrm{d}^{#1}{#2}/\mrm{d}{#3}^{#1}}} 
\newcommand{\lbrb}[1]{\ensuremath{\mkern-3mu\left[{#1}\right]\mkern-3mu}} 
\newcommand{\lbrrbr}[1]{\ensuremath{\mkern-3mu\left\{{#1}\right\}\mkern-3mu}} 
\newcommand{\lprp}[1]{\ensuremath{\mkern-3mu\left({#1}\right)\mkern-3mu}} 
\newcommand{\lprpsup}[2]{\ensuremath{\mkern-3mu\left({#1}\right)^{#2}\mkern-3mu}} 
\newcommand{\lbrbsup}[2]{\ensuremath{\mkern-3mu\left[{#1}\right]^{#2}\mkern-3mu}} 
\newcommand{\mad}[1]{\ensuremath{\genericmathaccenti{\acute}{#1}}} 
\newcommand{\madd}[1]{\ensuremath{\genericmathaccentii{\acute}{#1}}} 
\newcommand{\mb}[1]{\ensuremath{\bar{#1}}} 
\newcommand{\mbf}[1]{\ensuremath{\mathbf{#1}}} 
\newcommand{\mbfb}[1]{\ensuremath{\mb{\mathbf{#1}}}} 
\newcommand{\md}[1]{\ensuremath{\dot{#1}}} 
\newcommand{\mdd}[1]{\ensuremath{\ddot{#1}}} 
\newcommand{\mh}[1]{\ensuremath{\hat{#1}}} 
\newcommand{\mfkb}[1]{\ensuremath{\mb{\mfk{#1}}}} 
\newcommand{\pdrv}[2]{\ensuremath{\dfrac{\partial{#1}}{\partial{#2}}}} 
\newcommand{\mbb}[1]{\ensuremath{\mathbb{#1}}} 
\newcommand{\mcl}[1]{\ensuremath{\mathcal{#1}}} 
\newcommand{\mfk}[1]{\ensuremath{\mathfrak{#1}}} 
\newcommand{\mrm}[1]{\ensuremath{\mathrm{#1}}} 
\newcommand{\tbf}[1]{\textbf{#1}} 
\newcommand{\tc}{\text{c}} 
\newcommand{\tJ}{\text{J}} 
\newcommand{\tL}{\text{L}} 
\newcommand{\ts}{\text{s}} 
\newcommand{\tsp}{\text{sp}} 
\newcommand{\tB}{\text{B}} 
\newcommand{\mcolc}[2]{\multicolumn{#1}{c}{#2}} 
\newcommand{\mrow}[2]{\multirow{#1}{*}{#2}} 
\setlength{\arraycolsep}{1.5pt}
\newcommand{\MATLAB}{\texttt{MATLAB}} 
\newcommand{\ode}[1]{\texttt{ode{#1}}} 
\newcommand{\GPOPSII}{\mbb{GPOPS}\text{-}\mbb{II}} 
\newcommand{\enquote}[1]{``{#1},''} 

\begin{abstract}
The problem of minimizing the transfer time between periodic orbits in the Earth-Moon elliptic restricted three-body problem using a multi-mode propulsion system is considered.  By employing the true anomaly on the primary orbit as the independent variable and introducing normalized time as an additional state, the need to repeatedly solve Kepler's equation at arbitrary epochs is eliminated.  Furthermore, a propellant constraint is imposed on the high-thrust mode to activate the multi-mode capabilities of the system and balance efficiency with maneuverability.  The minimum-time optimal control problem is formulated as a three-phase trajectory consisting of a coast along the initial periodic orbit, a controlled transfer, and a coast along the terminal periodic orbit.  The three-phase optimal control problem is then solved using an adaptive Gaussian quadrature direct collocation method.  Case studies are presented for transfers from an $\tL_2$ southern halo orbit to a near-rectilinear halo orbit, analyzing the impact of different single- and multi-mode propulsion architectures and varying propellant constraint values. Finally, the methodology developed in this paper provides a systematic framework for generating periodic orbit transfers in three-body systems using single- and multi-mode propulsion systems.
\end{abstract}

\renewcommand{\baselinestretch}{1.5}
\normalsize\normalfont

\section{Introduction}\label{sec:intro}

The utilization of multi-mode propulsion systems onboard a spacecraft represents a frontier in astrodynamics and mission design.
Multi-mode propulsion systems integrate two or more propulsive modes with a shared propellant into a single system~\cite{RoveyBerg2020}, which allows the spacecraft to exploit the complementary advantages of distinct propulsion technologies within a single system.
Typical multi-mode systems combine a high-thrust, low specific impulse mode (for example, a chemical mode) for orbit maneuvers or large impulses with a low-thrust, high specific impulse mode (for example, an electric mode) for station-keeping.
As a result, multi-mode systems are of great interest due to their increased versatility and robustness to mission design requirements relative to a single-mode counterpart.
As space missions demand increasing payload capacity, reduced propellant consumption, and adaptable transfer strategies, multi-mode propulsion systems provide a critical pathway toward more efficient and resilient exploration architectures.
Thus, this research focuses on the use of multi-mode propulsion to accomplish periodic orbit transfers in three-body systems.

Within the realm of trajectory optimization, optimal control problems involving multi-mode propulsion systems are often solved using indirect shooting~\cite{TaheriGirard2020a,TaheriGirard2020b,AryaJunkins2021,ClineRovey2024b}.
For instance, Refs.~\cite{TaheriGirard2020a,TaheriGirard2020b} develop an indirect shooting framework to design optimal control strategies for problems with multiple modes of operation, where two-body and perturbed two-body dynamic models are employed to obtain optimal trajectories from Earth to Dionysus and the comet 67P.
Ref.~\cite{AryaJunkins2021} incorporates realistic multi-mode electric propulsion systems to demonstrate fuel-optimal gravity-assist trajectory optimization problem from Earth to the asteroid Psyche via a Mars flyby.
Similarly, Ref.~\cite{ClineRovey2024b} utilizes indirect shooting to obtain optimal propellant-constrained Earth to Mars and geosynchronous transfer orbit (GTO) to geostationary orbit (GEO) transfers using Modified Equinoctial Elements (MEEs) with perturbations.
In contrast, Ref.~\cite{ClineRovey2024a} utilizes a direct method within NASA's General Mission Analysis Tool (GMAT) to compute trajectories for lunar small satellite missions, including transfers to periodic orbits in the Earth–Moon system.
Outside the realm of trajectory optimization, multi-mode propulsion systems have been investigated for collision avoidance for spacecraft formations~\cite{Wang2017}.

While indirect shooting techniques remain the dominant approach for solving optimal orbit transfers using multi-mode propulsion systems, this work instead develops a direct collocation framework for optimal periodic orbit transfers.
Periodic orbits have supported long-duration scientific investigations for more than half a century~\cite{FarquharRichardson1977,SweetserWoodard2011,CheethamClarkson2022}, with renewed interest in recent years driven by cislunar exploration~\cite{MenzelBowers2023}.
Furthermore, direct collocation methods have been extensively applied within three-body systems to compute periodic orbits and design transfers between them in the circular restricted three-body problem (CR3BP)~\cite{OzimekHowell2009,OzimekHowell2010,GrebowHowell2011,HamanRao2024a,HamanRao2024b,HiraiwaBando2025}.
However, their extension to the elliptic restricted three-body problem (ER3BP) has received comparatively little attention, largely due to the additional complexity introduced by the explicit dependence on the independent variable in the ER3BP equations of motion (EOMs)~\cite{MoranteSoler2021}.
For instance, Ref.~\cite{OzimekHowell2009} employs odd-degree Gauss–Lobatto collocation schemes combined with Newton’s method to identify orbits enabling continuous surveillance of the lunar south pole using two-sided solar sail propulsion in the Earth-Moon CR3BP.
The study then demonstrates that quasi-periodic orbits, obtained under varying solar sail accelerations and regional constraints, can be successfully transitioned to ephemeris models, highlighting the robustness of direct methods with respect to naive initial guesses.
Similarly, Ref.~\cite{OzimekHowell2010} applies an adaptive odd-degree Gauss-Lobatto collocation method to design solar sail lunar pole-sitters in the Earth–Moon CR3BP using a global Fourier-series control law with constrained elevation angle and altitude, with the resulting periodic solutions verified using the \MATLAB{} ordinary differential equation (ODE) solvers \ode{45} and \ode{113}.
Ref.~\cite{GrebowHowell2011} further adapts the collocation approach from Ref.~\cite{OzimekHowell2009} to maximize pole-sitting coverage time when analyzing a three-phase low-thrust mission scenario.
Beyond lunar pole-sitters, Refs.~\cite{HamanRao2024a,HamanRao2024b} investigate minimum-time transfers between Earth–Moon $\tL_1$ and $\tL_2$ Lyapunov and halo orbits using a single-mode propulsion system, while Ref.~\cite{HiraiwaBando2025} employs time regularization and successive convexification techniques to construct minimum-time, low-thrust transfers between quasi-periodic orbits near $\tL_1$ and $\tL_2$.
In addition, Ref.~\cite{HiraiwaBando2025} neglects the spacecraft's mass variation due to propellant expenditure.
In contrast to these works, Ref.~\cite{DuLiu2023} uses indirect collocation techniques to compute low-energy transfers between periodic solutions in the Earth–Moon ER3BP.

This research is inspired by the work of Ref.~\cite{ClineRovey2024b}.
Although both Ref.~\cite{ClineRovey2024b} and this research focus on the design of optimal orbit transfers using multi-mode propulsion systems, the work presented in this paper significantly differs from Ref.~\cite{ClineRovey2024b} in several key aspects.
First, while Ref.~\cite{ClineRovey2024b} models the spacecraft's dynamics using perturbed two-body dynamics with eclipsing, this work employs the restricted three-body problem (R3BP) EOMs.
Second, integral constraints are introduced in this work to enforce an inequality constraint on the propellant consumed by the high-thrust mode, whereas Ref.~\cite{ClineRovey2024b} augments the cost functional with a penalty function term.
In addition, Ref.~\cite{ClineRovey2024b} employ a separate mass state for each propulsive mode, while this work uses a unified mass state for the spacecraft.
Third, the transfer scenarios differ substantially.
Ref.~\cite{ClineRovey2024b} examines single-phase GTO to GEO transfers; on the other hand, this work considers Earth-Moon R3BP periodic orbit transfers in a three-phase formulation.
In terms of methodology, Ref.~\cite{ClineRovey2024b} solves the optimal control problem using indirect shooting in combination with the \MATLAB{} \texttt{fsolve} function and ODE solver \ode{89}.
In contrast, all optimal trajectories in this work are obtained via an adaptive Legendre-Gauss Radau (LGR) direct collocation method~\cite{GargHuntington2010,GargRao2011b,PattersonRao2012}, which is chosen for its exponential convergence properties~\cite{HagerRao2016,HagerRao2017,HagerWang2018,HagerWang2019} and ability to handle highly nonlinear dynamics and a large number of constraints.
This approach eliminates the need to derive complex first-order optimality conditions or supply initial guesses for non-intuitive costate variables, as required by indirect methods.
Then, the resulting nonlinear programming problem (NLP) is solved using established optimization software~\cite{GillSaunders2005}, and all optimal trajectories are verified using the \MATLAB{} ODE solver \ode{113}~\cite{HamanRao2025}.

In addition to the aforementioned distinctions from existing literature, this paper makes several further contributions.
In this work, a novel three-phase, minimum-time, multi-mode optimal control problem is formulated within the R3BP, which comprises of the following: (1) a coast phase along the initial orbit, (2) a controlled transfer phase from the initial orbit to the terminal orbit, and (3) a coast phase along the terminal orbit.
Then, the controlled transfer phase is further partitioned into domains in order to optimize the switches in the throttles of each mode.
By employing the true anomaly on the primary orbit as the independent variable and introducing normalized time as an additional state variable, this novel formulation eliminates the need to repeatedly solve Kepler’s equation at arbitrary epochs in order to evaluate the corresponding minimum-time objective functional and any associated constraints.
Next, a numerical trajectory optimization study is performed for propellant-constrained, minimum-time transfers from an $\tL_2$ southern halo orbit to a near-rectilinear halo orbit (NRHO) in the Earth-Moon system.
Finally, the study examines various single- and multi-mode propulsion systems across a wide range of propellant constraint values, highlighting key features of the resulting optimal trajectories and the corresponding throttle structures for each mode.

This paper is organized as follows.
First, Section~\ref{sec:note} presents the notation and conventions utilized in this work.
Next, Section~\ref{sec:prob} formulates the three-phase optimal control problem for a minimum-time periodic orbit transfer in the Earth-Moon R3BP using a multi-mode propulsion system; Section~\ref{sec:prob} includes the modeling assumptions, the phase descriptions, as well as the necessary objective functional, dynamic and static parameter constraints, boundary conditions, and integral and path constraints.
The physical parameters and numerical approach utilized in this investigation are shown in Sections~\ref{sec:phys} and~\ref{sec:num}, respectively.
Then, Section~\ref{sec:res} provides and discusses results to the optimal control problem of Section~\ref{sec:prob}.
Finally, Section~\ref{sec:conc} provides conclusions on this research.


\section{Notation and Conventions}\label{sec:note}

In this paper, the following vector–matrix notation and conventions are used.
First, all scalars appear in lowercase and italic (for example, $a$), and all vectors are in $\mbb{R}^3$ unless explicitly stated otherwise and appear in lowercase and bold (for example, $\mbf{a}$).
All dimensional quantities are accented with a bar (for example, $\mb{a}$ or $\mbfb{a}$), whereas all nondimensional quantities are unmarked (for example, $a$ or $\mbf{a}$).
Next, the following function notation and conventions are used in this paper.
Let $f(x):\mbb{R}\rightarrow\mbb{R}$ denotes a scalar function $f$ of the independent variable $x$, where $f(x_i)\in\mbb{R}$ or $\left.f(x)\right|_{x=x_i}\in\mbb{R}$ denote the function $f(x)$ evaluated at $x=x_i$.
Then, $\mbf{f}(x):\mbb{R}\rightarrow\mbb{R}^n$ denotes a vector function $\mbf{f}$ of the independent variable $x$, where $\mbf{f}(x_i)\in\mbb{R}^n$ or $\left.\mbf{f}(x)\right|_{x=x_i}\in\mbb{R}^n$ denotes the function $\mbf{f}(x)$ evaluated at $x=x_i$.
Let $\mbf{f}(\mbf{x}):\mbb{R}^n\rightarrow\mbb{R}^m$ denote a vector function that maps $\mbf{x}\in\mbb{R}^n$ to $\mbf{f}(\mbf{x})\in\mbb{R}^m$.
Then, $\mbf{f}(\mbf{x},t):\mbb{R}^n\times\mbb{R}\rightarrow\mbb{R}^m$ is a multivariate vector function that maps $\mbf{x}\in\mbb{R}^n$ and $t\in\mbb{R}$ to $\mbf{f}(\mbf{x},t)\in\mbb{R}^m$.
Lastly, the following general notation and conventions are used.
Let $\md{x}\equiv\ildrv{x}{\mb{t}}$ and $\mdd{x}\equiv\ilndrv{2}{x}{\mb{t}}$ denote the first and second derivatives, respectively, of a quantity (in this case, $x$) with respect to dimensional time $\mb{t}$.
Furthermore, let $\mad{x}\equiv\ildrv{x}{\nu}$ and $\madd{x}\equiv\ilndrv{2}{x}{\nu}$ denote the first and second derivatives, respectively, of a quantity (in this case, $x$) with respect to the true anomaly $\nu$ on the primary system orbit.
Next, let $x_0$ and $x_f$ denote the initial and terminal values, respectively, of a quantity (in this case, $x$).
Finally, let $x^{\{p\}}$ denote a quantity (in this case, $x$), variable, or function in phase $p$, and let $x^*$ denote the optimal value of a quantity (in this case, $x$) or variable.


\section{Problem Description}\label{sec:prob}

In this section, the assumptions, objective functional, equations of motion, boundary conditions, and constraints that define the three-phase, minimum-time, multi-mode periodic orbit transfer optimal control problem are developed.
First, Section~\ref{sec:prob-ass} presents the assumptions used to model the motion of the spacecraft and primaries.
Next, Section~\ref{sec:prob-desc} describes each of the three phases of the optimal control problem.
Section~\ref{subsec:prob-obj} formulates the minimum-time objective functional, while Section~\ref{subsec:prob-eom} presents the equations of motion for the spacecraft as well as additional dynamic constraints.
Then, Sections~\ref{subsec:prob-s} and~\ref{subsec:prob-bnd} present the static parameter constraints and boundary conditions, respectively.
Sections~\ref{subsec:prob-int} and~\ref{subsec:prob-pth} present the integral and path constraints, respectively.
Finally, Section~\ref{subsec:prob-ocp} provides a complete description of the optimal control problem.

\subsection{Modeling Assumptions}\label{sec:prob-ass}

In this work, only closed conics within the R3BP are considered to model the motion of the primaries, which results in the elliptic R3BP (ER3BP).
In the ER3BP, the primary system orbit is assumed to be an ellipse (that is, $0\leq{e}<1$, where $e$ is the primary system orbit eccentricity).
Note that the circular R3BP (CR3BP) (that is, for $e=0$) is a special case of the ER3BP.
Using the true anomaly $\nu$ on the primary system orbit as the independent variable, the ER3BP EOMs are expressed in terms of nondimensionalized quantities, which helps reduce numerical error by bringing values to similar orders of magnitude and fixes the locations of the primaries in a non-uniformly rotating reference frame~\cite{Szebehely1967}.
In order to obtain nondimensionalized quantities, units of mass, length, and time are defined, which are all functions of parameters describing the primaries and the primary system orbit.
First, let the system mass ratio be defined as $\mu\equiv\mb{\mu}_2/(\mb{\mu}_1+\mb{\mu}_2)\leq{0.5}$, where $\mb{\mu}_i$ is the gravitational parameter of primary $i$, $(i=1,2)$.
Then, the mass unit (MU) is the sum of the masses of the primaries and defined as $\mfkb{M}\equiv(\mb{\mu}_1+\mb{\mu}_2)/\mb{G}$, where $\mfkb{M}$ is the characteristic mass and $\mb{G}$ is the universal gravitational constant.
Next, the length unit (LU) is the instantaneous distance between the primaries and defined as $\mfkb{L}(\nu)\equiv\mb{a}(1-e^2)/(1+e\cos(\nu))$, where $\mfkb{L}$ is the characteristic length and $\mb{a}$ is the primary system orbit semi-major axis.
Finally, the time unit (TU) instantaneously scales the nondimensionalized universal gravitational constant to unity and is defined as $\mfkb{T}(\nu)\equiv((\mfkb{L}(\nu))^3/(\mb{G}\mfkb{M}))^{1/2}$, where $\mfkb{T}$ is the characteristic time.
Note that the characteristic length and time in the non-circular ER3BP depend explicitly on $\nu$; however, the characteristic length and time are constant in the CR3BP.
From the mass, length, and time units, all other unit quantities (for example, speed and force) then follow.


\subsection{Phase Descriptions}\label{sec:prob-desc}

In this work, the problem of minimizing the transfer time between periodic orbits in the R3BP is formulated in three phases: (1) a coast phase along the initial CR3BP orbit, (2) a controlled ER3BP transfer phase from the initial orbit to the terminal orbit, and (3) a coast phase along the terminal CR3BP orbit.
Because no closed-form, analytic solution exists for the R3BP, numerical methods must be employed to determine a spacecraft's state history from an initial or terminal condition; thus, a set of position and velocity vector components is often utilized to define a R3BP periodic orbit.
Furthermore, a single parameter cannot be utilized to describe the location of a spacecraft on a R3BP periodic orbit at a specified time (that is, unlike the true anomaly on a Keplerian orbit), which demonstrates the need for the initial and terminal coast phases that allow the spacecraft's departure and arrival locations on the initial and terminal orbits, respectively, to be optimized.
Then, the three phases are linked together in chronological order as follows.
Phase 1 starts at the condition defining the initial periodic orbit and ends at the departure location from the initial periodic orbit.
Phase 2 starts at the departure location from the initial periodic orbit and ends at the arrival location on the terminal periodic orbit.
Phase 3 starts at the arrival location on the terminal periodic orbit and ends at the condition defining the terminal periodic orbit.


\subsection{Objective Functional}\label{subsec:prob-obj}

In this section, the objective functional for the R3BP minimum-time, multi-mode problem formulation is presented, where the goal is to minimize the duration of the controlled transfer phase.
In order to eliminate the need to evaluate Kepler's equation, remove the additional associated computational expense, and formulate a general minimum-time objective functional in the R3BP, normalized time $\tau$ is introduced as an additional state component.
Then, dimensional time $\mb{t}$ is scaled via $\tau\equiv\mb{t}/\mfkb{T}(0)$ in order to bring the value to a similar order of magnitude to the instantaneously nondimensionalized quantities, where $\mfkb{T}(0)$ (that is, the characteristic time evaluated at $\nu=0$) is a constant.
Note that the scaling characteristic time can be evaluated using any $\nu\in[0,2\pi]$; however, a value of $\nu=0$ is chosen for convenience.
Then, the minimum-time objective functional is written as
\begin{equation}\label{eq:prob-obj-dur}
    \mcl{J}=\tau^{\{2\}}(\nu^{\{2\}}_f)-\tau^{\{2\}}(\nu^{\{2\}}_0),
\end{equation}
which is the difference between the initial and terminal normalized time in the transfer phase.
Thus, the objective functional of Eq.~\eqref{eq:prob-obj-dur} can be evaluated directly using only the corresponding values of normalized time and does not require any evaluation of Kepler's equation.


\subsection{Equations of Motion}\label{subsec:prob-eom}

In this section, the equations of motion (that is, the dynamic constraints) for the R3BP minimum-time, multi-mode problem formulation are presented.
First, the dynamic constraints of the controlled ER3BP transfer phase are addressed.
The ER3BP equations of motion are expressed in terms of nondimensionalized pulsating coordinates in a non-uniformly rotating reference frame using the true anomaly as the independent variable~\cite{Szebehely1967}.
Thus, the differential equations of motion for the spacecraft under the gravitational influence of the two primaries and a perturbing thrust acceleration (that is, the controlled ER3BP EOMs) are given by
\begin{align}
	\label{eq:prob-eom-x}
	\madd{x}-2\mad{y} & =\grad{x}{\mfk{W}}+a_T\mh{u}_x, \\
	\label{eq:prob-eom-y}
	\madd{y}+2\mad{x} & =\grad{y}{\mfk{W}}+a_T\mh{u}_y, \\
	\label{eq:prob-eom-z}
	\madd{z}          & =\grad{z}{\mfk{W}}+a_T\mh{u}_z,
\end{align}
where $\{x,y,z\}$, $\{\mad{x},\mad{y},\mad{z}\}$, and $\{\madd{x},\madd{y},\madd{z}\}$ are the sets of the spacecraft's position, velocity, and acceleration components, respectively, and $\{\mh{u}_x,\mh{u}_y,\mh{u}_z\}$ is the associated set of unit-vector control direction components.
For use in the controlled ER3BP EOMs, the ER3BP pseudopotential function is defined as
\begin{equation}\label{eq:prob-eom-W}
	\mfk{W}(x,y,z,\nu)\equiv\dfrac{1}{1+e\cos(\nu)}\lbrb{\dfrac{1}{2}(x^2+y^2-e\cos(\nu)z^2)+\dfrac{1-\mu}{r_1}+\dfrac{\mu}{r_2}},
\end{equation}
the nondimensional distances between the two primaries and the spacecraft are
\begin{align}
	\label{eq:prob-eom-r1}
	r_1(x,y,z) & =\lbrbsup{(x+\mu)^2+y^2+z^2}{1/2}, \\
	\label{eq:prob-eom-r2}
	r_2(x,y,z) & =\lbrbsup{(x-1+\mu)^2+y^2+z^2}{1/2},
\end{align}
respectively, and the perturbing thrust acceleration magnitude is given by
\begin{equation}\label{eq:prob-eom-aT}
    a_T(m_{\ts},\delta_1,\delta_2,\nu)\equiv\dfrac{1}{1+e\cos(\nu)}\dfrac{1}{\mb{m}_{\ts,0}m_{\ts}}\dfrac{(\mfkb{T}(\nu))^2}{\mfkb{L}(\nu)}(\delta_1\mb{T}_{\max,1}+\delta_2\mb{T}_{\max,2}).
\end{equation}
In Eq.~\eqref{eq:prob-eom-aT}, $\delta_i$ is the throttle for mode $i$, $(i=1,2)$, and $\mb{T}_{\max,i}$ is the corresponding maximum thrust magnitude.
As discussed in Section~\ref{subsec:prob-obj}, normalized time is introduced as an additional state in the transfer phase; thus, the corresponding dynamic constraint is given by
\begin{equation}\label{eq:prob-eom-tau}
    \mad{\tau}=\dfrac{1}{\md{\nu}\mfkb{T}(0)},
\end{equation}
where the magnitude of the dimensional angular velocity of the rotating reference frame observed in the inertial reference frame is given by
\begin{equation}\label{eq:prob-eom-angvel}
    \md{\nu}(\nu)\equiv\mb{n}\dfrac{(1+e\cos(\nu))^2}{(1-e^2)^{3/2}},
\end{equation}
where $\mb{n}\equiv(\mb{G}\mfkb{M}/\mb{a}^3)^{1/2}$ is the dimensional mean motion of the primary system orbit.
To account for any variation in the spacecraft's mass due to any propellant expenditure, the spacecraft's normalized mass dynamics is given as 
\begin{equation}\label{eq:prob-eom-m}
    \mad{m}_{\ts}=-\dfrac{1}{\md{\nu}\mb{g}_0\mb{m}_{\ts,0}}\lprp{\dfrac{\mb{T}_{\max,1}}{\mb{I}_{\tsp,1}}\delta_1+\dfrac{\mb{T}_{\max,2}}{\mb{I}_{\tsp,2}}\delta_2},
\end{equation}
where the spacecraft's dimensional mass $\mb{m}_{\ts}$ is normalized via $m_{\ts}=\mb{m}_{\ts}/\mb{m}_{\ts,0}\in(0,1]$, $\mb{m}_{\ts,0}$ is the spacecraft's initial mass, $\mb{I}_{\tsp,i}$ is the specific impulse for mode $i$, $(i=1,2)$, and $\mb{g}_0$ is the acceleration due to gravity at sea level.
It is important to note that dynamic constraints of Eqs.~\eqref{eq:prob-eom-x}--\eqref{eq:prob-eom-z},~\eqref{eq:prob-eom-tau}, and~\eqref{eq:prob-eom-m} do not require any evaluation of Kepler's equation.
In the controlled transfer phase (that is, in phase $p=2$), the variables $\mbf{x}^{\{2\}}\equiv[x,y,z,\mad{x},\mad{y},\mad{z},\tau,m_{\ts}]$ and $\mbf{u}^{\{2\}}\equiv[\mh{u}_x,\mh{u}_y,\mh{u}_z,\delta_1,\delta_2]$ comprise the state and control, respectively.
Finally, the dynamic constraints of the initial and terminal CR3BP coast phases are addressed.
The ballistic CR3BP EOMs are obtained by setting $e=0$ and $a_T=0$ in Eqs.~\eqref{eq:prob-eom-x}--\eqref{eq:prob-eom-z} and~\eqref{eq:prob-eom-tau}, where \textit{Jacobi's constant} (that is, a CR3BP integral of motion) is defined as $\mfk{J}\equiv\left.2\mfk{W}\right|_{e=0}-(\mad{x}^2+\mad{y}^2+\mad{z}^2)$.
Then, the variables $\mbf{x}^{\{p\}}\equiv[x,y,z,\mad{x},\mad{y},\mad{z},\tau]$, $(p=1,3)$, comprise the state in the coast phases.
In order to simplify the problem formulation and reduce the size of the resulting NLP, control and the spacecraft's normalized mass are not included in the coast phases.


\subsection{Static Parameter Constraints}\label{subsec:prob-s}

In this section, the static parameter constraints for the R3BP minimum-time, multi-mode problem formulation are presented.
First, unique departure and arrival locations on the initial and terminal CR3BP periodic orbits, respectively, must be determined.
Similar to the true anomaly on a Keplerian orbit, the set of static parameters $\{\tc_0,\ts_0\}$ determines the duration of the coast phase along the initial orbit and, thus, the departure location on the initial orbit, where $\tc_0$ and $\ts_0$ can be interpreted as the cosine and sine, respectively, of some angle defining the departure location on the initial orbit.
A unique coast duration along the initial orbit is then ensured by the constraint
\begin{equation}\label{eq:prob-s-cs0}
    \tc^2_0+\ts^2_0-1=0.
\end{equation}
Then, a unique coast duration along the terminal orbit is ensured using the set of static parameters $\{\tc_f,\ts_f\}$ by the constraint
\begin{equation}\label{eq:prob-s-csf}
    \tc^2_f+\ts^2_f-1=0.
\end{equation}
In a similar manner, a unique initial true anomaly on the primary system orbit is ensured using the set of static parameters $\{\tc_{\nu_0},\ts_{\nu_0}\}$ by the constraint
\begin{equation}\label{eq:prob-s-csv0}
    \tc^2_{\nu_0}+\ts^2_{\nu_0}-1=0.
\end{equation}
The constraints of Eqs.~\eqref{eq:prob-s-cs0}--\eqref{eq:prob-s-csv0} also ensure $\mbf{s}\equiv\{\tc_0,\ts_0,\tc_f,\ts_f,\tc_{\nu_0},\ts_{\nu_0}\}\in[-1,+1]$.
Finally, all static parameters are free variables in the problem and further discussed in Section~\ref{subsec:prob-bnd}.


\subsection{Boundary Conditions}\label{subsec:prob-bnd}

In this section, the boundary conditions for the R3BP minimum-time, multi-mode problem formulation are presented.
First, the initial and terminal coast phase boundary conditions are addressed.
The spacecraft must depart from the desired initial orbit, which is enforced by the constraint
\begin{equation}\label{eq:prob-bnd-x1_0}
    \mbf{x}^{\{1\}}(\nu^{\{1\}}_0)=[x_0,y_0,z_0,\mad{x}_0,\mad{y}_0,\mad{z}_0,\tau_0]=\mbf{x}_0.
\end{equation}
In Eq.~\eqref{eq:prob-bnd-x1_0}, the variables $\{x_0,y_0,z_0,\mad{x}_0,\mad{y}_0,\mad{z}_0\}$ are fixed and determined by the specified initial periodic orbit.
Because $\tau$ does not appear explicitly in the EOMs of Section~\ref{subsec:prob-eom}, the initial value does not impact the solution; however, the initial value is fixed at $\tau_0=0$ for convenience.
Then, the initial coast duration is bounded to the interval $[0,\mfkb{p}_0]$ via
\begin{equation}\label{eq:prob-bnd-pdur0}
    \tau^{\{1\}}(\nu^{\{1\}}_f)-\tau^{\{1\}}(\nu^{\{1\}}_0)=\lprp{\dfrac{\atantwo(\ts_0,\tc_0)+\pi}{2\pi}}\dfrac{\mfkb{p}_0}{\mfkb{T}(0)},
\end{equation}
where $\mfkb{p}_0$ is the period of the initial orbit.
Note that the four-quadrant inverse tangent function $\atantwo$ returns a value on the interval $[-\pi,\pi]$, which is then scaled to the interval $[0,1]$ in Eq.~\eqref{eq:prob-bnd-pdur0}.
As a result, the initial coast phase traces the initial orbit from the initial state defined by Eq.~\eqref{eq:prob-bnd-x1_0} for a duration specified by Eq.~\eqref{eq:prob-bnd-pdur0} until the departure location on the initial orbit (that is, when the spacecraft enters the transfer phase).
Similarly, the spacecraft must arrive on the desired terminal orbit, which is enforced by the constraint
\begin{equation}\label{eq:prob-bnd-x3_f}
    \mbf{x}^{\{3\}}(\nu^{\{3\}}_f)=[x_f,y_f,z_f,\mad{x}_f,\mad{y}_f,\mad{z}_f,\tau_f]=\mbf{x}_3.
\end{equation}
In Eq.~\eqref{eq:prob-bnd-x3_f}, the variables $\{x_f,y_f,z_f,\mad{x}_f,\mad{y}_f,\mad{z}_f\}$ are fixed and determined by the specified terminal periodic orbit, and $\tau_f$ is free.
Then, the terminal coast duration is bounded to the interval $[0,\mfkb{p}_f]$ via
\begin{equation}\label{eq:prob-bnd-pdur3}
    \tau^{\{3\}}(\nu^{\{3\}}_f)-\tau^{\{3\}}(\nu^{\{3\}}_0)=\lprp{\dfrac{\atantwo(\ts_f,\tc_f)+\pi}{2\pi}}\dfrac{\mfkb{p}_f}{\mfkb{T}(0)},
\end{equation}
where $\mfkb{p}_f$ is the period of the terminal orbit.
As a result, the terminal coast phase traces the terminal orbit from the arrival location on the terminal orbit (that is, when the spacecraft leaves the transfer phase) for a duration specified by Eq.~\eqref{eq:prob-bnd-pdur3} until the terminal state defined by Eq.~\eqref{eq:prob-bnd-x3_f}.

Next, continuity across phases is addressed along with the boundary conditions in the transfer phase.
Because the three phases occur in chronological order, continuity of the independent variable across phases is enforced by the constraint
\begin{equation}\label{eq:prob-bnd-nucnt123}
    \nu^{\{p\}}_f=\nu^{\{p+1\}}_0=\nu_p, \qquad (p=1,2).
\end{equation}
Note that $\nu_p$, $(p=1,2)$, are free variables in the problem and correspond to the true anomaly on the primary system orbit at the interface between phases.
Because $\nu$ appears explicitly in the EOMs of Section~\ref{subsec:prob-eom}, the initial true anomaly is determined by the corresponding static parameters via
\begin{equation}\label{eq:prob-bnd-nu0}
    \nu^{\{1\}}_0=\atantwo(\ts_{\nu_0},\tc_{\nu_0})+\pi=\nu_0,
\end{equation}
where  $\nu_0\in[0,2\pi]$.
Equation~\eqref{eq:prob-bnd-nu0} also determines the locations of the primaries at the start of the transfer phase, which can significantly impact the solution when $0<e<1$.
Furthermore, the terminal true anomaly $\nu^{\{3\}}_f=\nu_3$ is free.
In a similar manner, continuity of the state across phases is also enforced.
In order to ensure that the spacecraft's dimensional position vector is continuous across phases, position continuity from the terminus of the initial CR3BP coast phase to the start of the ER3BP transfer phase requires
\begin{equation}\label{eq:prob-bnd-xyzcnt12}
    \lbrb{x^{\{1\}}(\nu^{\{1\}}_f),y^{\{1\}}(\nu^{\{1\}}_f),z^{\{1\}}(\nu^{\{1\}}_f)}=\gamma(\nu^{\{2\}}_0)\lbrb{x^{\{2\}}(\nu^{\{2\}}_0),y^{\{2\}}(\nu^{\{2\}}_0),z^{\{2\}}(\nu^{\{2\}}_0)},
\end{equation}
and position continuity from the terminus of the ER3BP transfer phase to the start of the terminal CR3BP coast phase requires
\begin{equation}\label{eq:prob-bnd-xyzcnt23}
    \lbrb{x^{\{3\}}(\nu^{\{3\}}_0),y^{\{3\}}(\nu^{\{3\}}_0),z^{\{3\}}(\nu^{\{3\}}_0)}=\gamma(\nu^{\{2\}}_f)\lbrb{x^{\{2\}}(\nu^{\{2\}}_f),y^{\{2\}}(\nu^{\{2\}}_f),z^{\{2\}}(\nu^{\{2\}}_f)},
\end{equation}
where $\gamma(\nu)\equiv(1-e^2)/(1+e\cos(\nu))$ accounts for the different CR3BP and ER3BP characteristic lengths.
Similarly, the basic kinematic equation is employed to ensure that the spacecraft's dimensional velocity vectors are continuous across phases.
As a result, velocity continuity from the terminus of the initial CR3BP coast phase to the start of the ER3BP transfer phase requires
\begin{equation}\label{eq:prob-bnd-vxvyvzcnt12}
    \dfrac{\mb{a}\mb{n}}{\mfkb{L}(\nu^{\{2\}}_0)\md{\nu}(\nu^{\{2\}}_0)}\lbrb{%
    \begin{array}{c}
        \mad{x}^{\{1\}}(\nu^{\{1\}}_f) \\
        \mad{y}^{\{1\}}(\nu^{\{1\}}_f) \\
        \mad{z}^{\{1\}}(\nu^{\{1\}}_f) \\
    \end{array}
    }=\left.\lbrrbr{\lbrb{%
    \begin{array}{c}
        \mad{x}^{\{2\}}(\nu) \\
        \mad{y}^{\{2\}}(\nu) \\
        \mad{z}^{\{2\}}(\nu) \\
    \end{array}
    }+\eta(\nu)\lbrb{%
    \begin{array}{c}
        x^{\{2\}}(\nu) \\
        y^{\{2\}}(\nu) \\
        z^{\{2\}}(\nu) \\
    \end{array}
    }+\xi(\nu)\lbrb{%
    \begin{array}{c}
        -y^{\{2\}}(\nu) \\
        x^{\{2\}}(\nu) \\
        0 \\
    \end{array}
    }}\right|_{\nu=\nu^{\{2\}}_0},\!\!
\end{equation}
and velocity continuity from the terminus of the ER3BP transfer phase to the start of the terminal CR3BP coast phase requires
\begin{equation}\label{eq:prob-bnd-vxvyvzcnt23}
    \dfrac{\mb{a}\mb{n}}{\mfkb{L}(\nu^{\{2\}}_f)\md{\nu}(\nu^{\{2\}}_f)}\lbrb{%
    \begin{array}{c}
        \mad{x}^{\{3\}}(\nu^{\{3\}}_0) \\
        \mad{y}^{\{3\}}(\nu^{\{3\}}_0) \\
        \mad{z}^{\{3\}}(\nu^{\{3\}}_0) \\
    \end{array}
    }=\left.\lbrrbr{\lbrb{%
    \begin{array}{c}
        \mad{x}^{\{2\}}(\nu) \\
        \mad{y}^{\{2\}}(\nu) \\
        \mad{z}^{\{2\}}(\nu) \\
    \end{array}
    }+\eta(\nu)\lbrb{%
    \begin{array}{c}
        x^{\{2\}}(\nu) \\
        y^{\{2\}}(\nu) \\
        z^{\{2\}}(\nu) \\
    \end{array}
    }+\xi(\nu)\lbrb{%
    \begin{array}{c}
        -y^{\{2\}}(\nu) \\
        x^{\{2\}}(\nu) \\
        0 \\
    \end{array}
    }}\right|_{\nu=\nu^{\{2\}}_f},\!\!
\end{equation}
where $\eta(\nu)\equiv{}e\sin(\nu)/(1+e\cos(\nu))$ and $\xi(\nu)\equiv{}1-\mb{n}/\md{\nu}(\nu)$.
Note that Eqs.~\eqref{eq:prob-bnd-vxvyvzcnt12} and~\eqref{eq:prob-bnd-vxvyvzcnt23} account for the different CR3BP and ER3BP angular velocity magnitudes shown in Eq.~\eqref{eq:prob-eom-angvel}.
Next, continuity of normalized time across phases is enforced by the constraint
\begin{equation}\label{eq:prob-bnd-tauucnt123}
    \tau^{\{p\}}(\nu^{\{p\}}_f)=\tau^{\{p+1\}}(\nu^{\{p+1\}}_0)=\tau_p, \qquad (p=1,2).
\end{equation}
where $\tau_p$, $(p=1,2)$, are free variables in the problem and correspond to the normalized time at the interface between phases.
It is important to note that all position, velocity, and normalized time variables at the interface between phases are free.
Although ER3BP periodic orbits are not considered in this work, such orbits can be integrated in a straightforward manner in this problem formulation by modifying the boundary conditions of Eqs.~\eqref{eq:prob-bnd-xyzcnt12}--\eqref{eq:prob-bnd-vxvyvzcnt23} and updating the associated dynamic constraints.
Finally, at the start of the transfer phase, the spacecraft's dimensional mass is $\mb{m}_{\ts,0}$; thus, the spacecraft's normalized mass at the start of the transfer phase is
\begin{equation}\label{eq:prob-bnd-m0}
    m^{\{2\}}_{\ts}(\nu^{\{2\}}_0)=1,
\end{equation}
and $m^{\{2\}}_{\ts}(\nu^{\{2\}}_f)$ is free.


\subsection{Integral Constraints}\label{subsec:prob-int}

In this section, the integral constraints for the R3BP minimum-time problem, multi-mode formulation are presented.
A propellant constraint on mode 1 (that is, the high-thrust mode) excites the multi-mode propulsion capabilities and introduces switches in the throttle structure of each mode~\cite{ClineRovey2024b}.
In this work, integral constraints are utilized to enforce an inequality constraint on the propellant consumed by mode 1, which is different from Ref.~\cite{ClineRovey2024b} that augment the cost functional with a penalty function term.
Because the spacecraft's normalized mass flow rate is a function of both modes, a lower bound on the spacecraft's normalized mass does not directly limit the amount of propellant consumed by mode 1 only.
In order to enforce an upper bound on the allowable propellant consumption of mode 1, the amount of propellant consumed by mode 1 must first be determined.
Utilizing the spacecraft's normalized mass dynamics of Eq.~\eqref{eq:prob-eom-m} and ignoring the contribution by mode 2, the integral in the transfer phase is obtained via
\begin{equation}\label{eq:prob-int-m1}
    q^{\{2\}}_{m_{\ts{1}}}=\int_{\nu^{\{2\}}_0}^{\nu^{\{2\}}_f}\lprp{\dfrac{1}{\md{\nu}\mb{g}_0\mb{m}_{\ts,0}}\dfrac{\mb{T}_{\max,1}}{\mb{I}_{\tsp,1}}\delta_1}\dnu,
\end{equation}
where $q^{\{2\}}_{m_{\ts{1}}}$ is equivalent to the amount of propellant consumed by mode 1 in the transfer phase.
Then, an upper bound on the total propellant consumption of mode 1 is enforced via
\begin{equation}\label{eq:prob-int-m1sum}
    0\leq{q^{\{2\}}_{m_{\ts{1}}}}\leq{m_{\ts{1},u}},
\end{equation}
where $m_{\ts{1},u}$ denotes the upper bound on the propellant consumption of mode 1.
It is important to note that Eqs.~\eqref{eq:prob-int-m1} and~\eqref{eq:prob-int-m1sum} are employed using the spacecraft's normalized mass (for example, a value of $m_{\ts{1},u}=0.1$ corresponds to a 10 [kg] propellant constraint for a 100 [kg] spacecraft).


\subsection{Path Constraints}\label{subsec:prob-pth}

In this section, the path constraints for the R3BP minimum-time, multi-mode problem formulation are presented, which are all only enforced in the controlled transfer phase.
First, a unit-vector thrust direction requires
\begin{equation}\label{eq:prob-pth-unit}
    \mh{u}_x^2+\mh{u}_y^2+\mh{u}_z^2-1=0,
\end{equation}
which also ensures $\{\mh{u}_x,\mh{u}_y,\mh{u}_z\}\in[-1,+1]$.
Similar to Ref.~\cite{ClineRovey2024b}, the constraint enabling multi-mode propulsion capabilities is defined as
\begin{equation}\label{eq:prob-pth-one}
    \delta_1\delta_2=0,
\end{equation}
which ensures at least one of the two thrusting modes is inactive at all times (that is, at most one mode can be active at all times).
Then, the throttle for each mode is bounded via
\begin{equation}\label{eq:prob-pth-ubnd}
    0\leq\delta_i\leq{1}, \qquad (i=1,2),
\end{equation}
which ensures that the resultant thrust magnitude from mode $i$ is bounded between zero and $\mb{T}_{\max,i}$, $(i=1,2)$.
Numerical difficulties often arise in close proximity to the primaries due to the existence of singularities in the dynamic model at the locations of the primaries (that is, at $x=-\mu$ and $x=1-\mu$).
Finally, to help avoid such numerical difficulties as well as ensure that the spacecraft does not collide with the primaries, minimum altitude constraints are enforced via
\begin{equation}\label{eq:prob-pth-alt}
    \mb{R}_i+\mb{h}_{i,\min}-\mfkb{L}(\nu)r_i\leq{0}, \qquad (i=1,2),
\end{equation}
where $r_i$ is the distance between the spacecraft and primary $i$ given by either Eq.~\eqref{eq:prob-eom-r1} or~\eqref{eq:prob-eom-r2}, and $\mb{R}_i$ and $\mb{h}_{i,\min}$ are the dimensional mean equatorial radius and specified minimum altitude corresponding to primary $i$, $(i=1,2)$, respectively.


\subsection{Optimal Control Problem}\label{subsec:prob-ocp}

The R3BP minimum-time, multi-mode periodic orbit transfer optimal control problem is formally stated as follows.
Determine the state $\mbf{x}^{\{p\}}\in\mbb{R}^7$, $(p=1,3)$, in the initial and terminal coast phases, the state $\mbf{x}^{\{2\}}\in\mbb{R}^8$, control $\mbf{u}^{\{2\}}\in\mbb{R}^5$, and integral $q^{\{2\}}_{m_{\ts{1}}}\in\mbb{R}$ in the transfer phase, the initial true anomalies $\nu^{\{p\}}_0\in\mbb{R}$ and terminal true anomalies $\nu^{\{p\}}_f\in\mbb{R}$ in all phases $p$, $(p=1,2,3)$, as well as the static parameters $\mbf{s}\in\mbb{R}^6$ that minimize the objective functional of Eq.~\eqref{eq:prob-obj-dur} in Section~\ref{subsec:prob-obj} subject to the dynamic constraints, static parameter constraints, boundary conditions, integral constraints, and path constraints of Sections of Section~\ref{subsec:prob-eom}--\ref{subsec:prob-pth}, respectively.



\section{Physical Parameters}\label{sec:phys}

In this section, the physical parameters utilized in the R3BP minimum-time, multi-mode optimal control problem formulation of Section~\ref{sec:prob} are presented.
First, the Earth-Moon R3BP model parameters are discussed in Section~\ref{subsec:phys-model}.
Then, the spacecraft parameters are shown in Section~\ref{subsec:phys-sc}.
Finally, Section~\ref{subsec:phys-orbit} presents the target periodic orbits.

\subsection{Model Parameters}\label{subsec:phys-model}

In this section, the physical model parameters utilized in the R3BP minimum-time, multi-mode optimal control problem formulation of Section~\ref{sec:prob} are presented.
For this investigation, the Earth-Moon-spacecraft three-body system is of interest, where the Earth-Moon system parameters are shown in Table~\ref{tab:phys-model-sys}.
Furthermore, minimum altitude values of $\mb{h}_{1,\min}=500$ [km] and $\mb{h}_{2,\min}=200$ [km] are utilized in the path constraint of Eq.~\eqref{eq:prob-pth-alt}.
\begin{table}[!t]
    \centering
    \caption{Earth-Moon system parameters.}
    \begin{tabular*}{0.60\textwidth}{@{\extracolsep\fill}crc@{}} \hline\hline
        \tbf{Parameter} & \mcolc{1}{\tbf{Value}} & \tbf{Units} \\ \hline
        $\mb{\mu}_1$ & 3.986002221116981$\times{10}^{+5~}$ & [km$^3$/s$^2$] \\
        $\mb{\mu}_2$ & 4.902797989481230$\times{10}^{+3~}$ & [km$^3$/s$^2$] \\
        $\mb{R}_1$   & 6.378137$\times{10}^{+3~}$          & [km] \\
        $\mb{R}_2$   & 1.7371$\times{10}^{+3~}$            & [km] \\
        $\mb{G}$     & 6.67430$\times{10}^{-20}$           & [km$^3$/kg/s$^2$] \\
        $\mb{g}_0$   & 9.80665$\times{10}^{-3~}$           & [km/s$^2$] \\
        $\mu$        & 1.215058560962404$\times{10}^{-2~}$ & [-] \\
        $\mb{a}$     & 3.89703$\times{10}^{+5~}$           & [km] \\
        $e$          & 5.49$\times{10}^{-2~}$              & [-] \\ \hline\hline
    \end{tabular*}
    \label{tab:phys-model-sys}
\end{table}


\subsection{Spacecraft Parameters}\label{subsec:phys-sc}

In this section, the physical spacecraft parameters utilized in the R3BP minimum-time, multi-mode optimal control problem formulation of Section~\ref{sec:prob} are presented.
For this investigation, the spacecraft and corresponding multi-mode propulsion system specifications are shown in Tables~\ref{tab:phys-sc-ts1} and~\ref{tab:phys-sc-ts2}, where the Table~\ref{tab:phys-sc-ts1} specifications are taken directly from Ref.~\cite{ClineRovey2024b}.
Mode 1 is assumed to be a 1 [N], 250 [s] mode, whereas mode 2 is assumed to be a 0.5 [N], 3100 [s] mode.
As mentioned in Ref.~\cite{ClineRovey2024b}, the physical parameters shown in Table~\ref{tab:phys-sc-ts1} are representative of a monopropellant thruster for mode 1 and a high-power thruster for mode 2, which requires 13 [kW] of power assuming 60\% efficiency.
Furthermore, Ref.~\cite{ClineRovey2024b} states that mode 2 is currently infeasible for a 100 [kg] spacecraft; however, the specifications are chosen to provide insight into the structure of the solution at the expense of some realism.
Finally, the specifications of Table~\ref{tab:phys-sc-ts2} are identical to that shown in Table~\ref{tab:phys-sc-ts1} except that $\mb{T}_{\max,2}$ is decreased from 0.5 [N] to 0.25 [N] to provide insight into the solution structure as the maximum thrust from mode 2 is modified.
\begin{table}[!t]
    \centering
    \begin{minipage}{0.48\textwidth}
        \centering
        \caption{Spacecraft parameters (Case 1).}
        \begin{tabular*}{0.8\textwidth}{@{\extracolsep\fill}crc@{}} \hline\hline
            \tbf{Parameter} & \mcolc{1}{\tbf{Value}} & \tbf{Units} \\ \hline
            $\mb{T}_{\max,1}$ & 1.0$\times{10}^{-3}$ & [kN] \\
            $\mb{I}_{\tsp,1}$ & 2.5$\times{10}^{+2}$ & [s] \\
            $\mb{T}_{\max,2}$ & 0.5$\times{10}^{-3}$ & [kN] \\
            $\mb{I}_{\tsp,2}$ & 3.1$\times{10}^{+3}$ & [s] \\
            $\mb{m}_{\ts,0}$  & 1.0$\times{10}^{+2}$ & [kg] \\ \hline\hline
        \end{tabular*}
        \label{tab:phys-sc-ts1}
    \end{minipage}
    \begin{minipage}{0.48\textwidth}
        \centering
        \caption{Spacecraft parameters (Case 2).}
        \begin{tabular*}{0.8\textwidth}{@{\extracolsep\fill}crc@{}} \hline\hline
            \tbf{Parameter} & \mcolc{1}{\tbf{Value}} & \tbf{Units} \\ \hline
            $\mb{T}_{\max,1}$ & 1.0$\times{10}^{-3}$ & [kN] \\
            $\mb{I}_{\tsp,1}$ & 2.5$\times{10}^{+2}$ & [s] \\
            $\mb{T}_{\max,2}$ & 2.5$\times{10}^{-4}$ & [kN] \\
            $\mb{I}_{\tsp,2}$ & 3.1$\times{10}^{+3}$ & [s] \\
            $\mb{m}_{\ts,0}$  & 1.0$\times{10}^{+2}$ & [kg] \\ \hline\hline
        \end{tabular*}
        \label{tab:phys-sc-ts2}
    \end{minipage}
\end{table}

\subsection{Orbit Parameters}\label{subsec:phys-orbit}

\begin{figure}[!t]
    \centering
    \begin{minipage}{0.48\textwidth}
        \centering
        \subfloat[3-dimensional.]{%
            \includegraphics[height=18.4em]{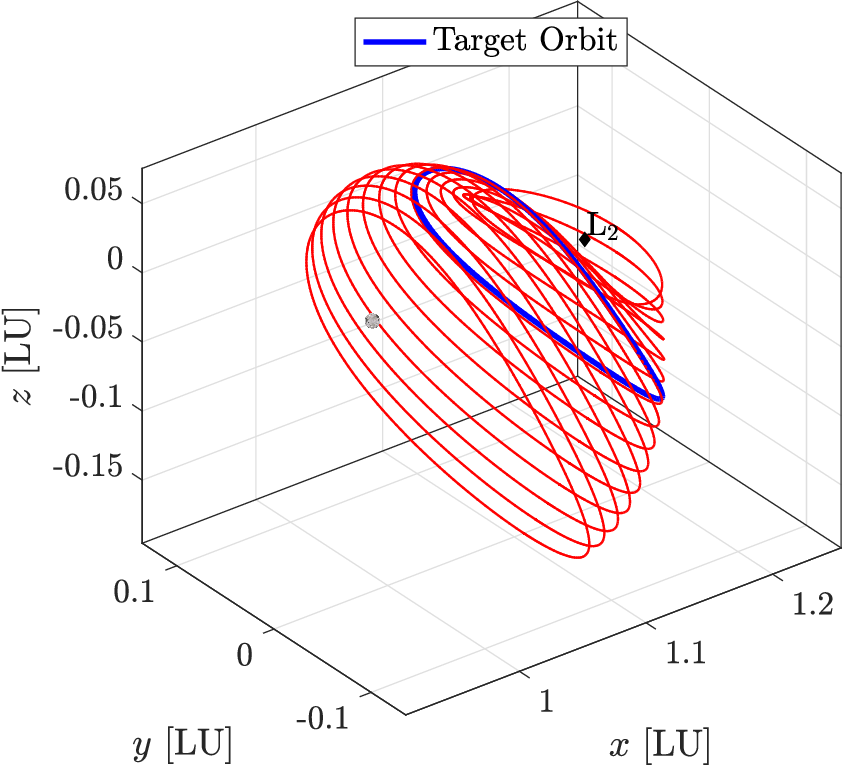}%
            \label{fig:phys-orbit-l2shalo_3d}
            }
        \\
        \hspace*{\fill}%
        \subfloat[$xy$-plane.]{%
            \includegraphics[height=5.1em]{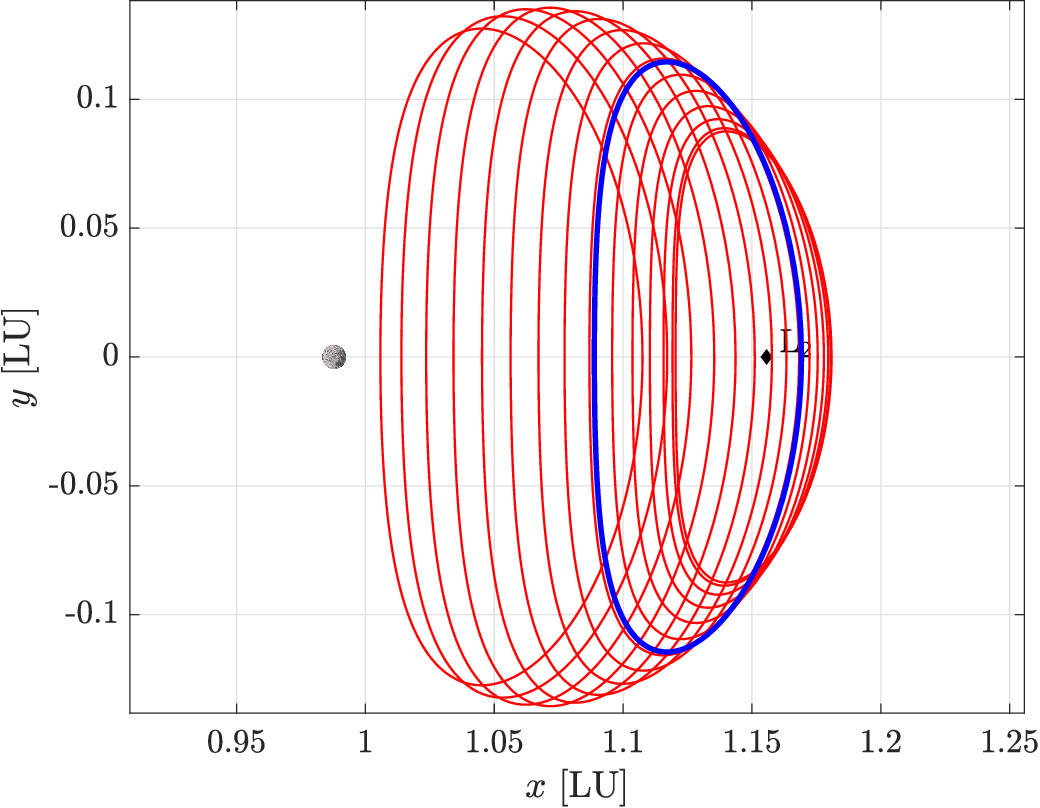}%
            \label{fig:phys-orbit-l2shalo_xy}
            }
        \subfloat[$xz$-plane.]{%
            \includegraphics[height=5.1em]{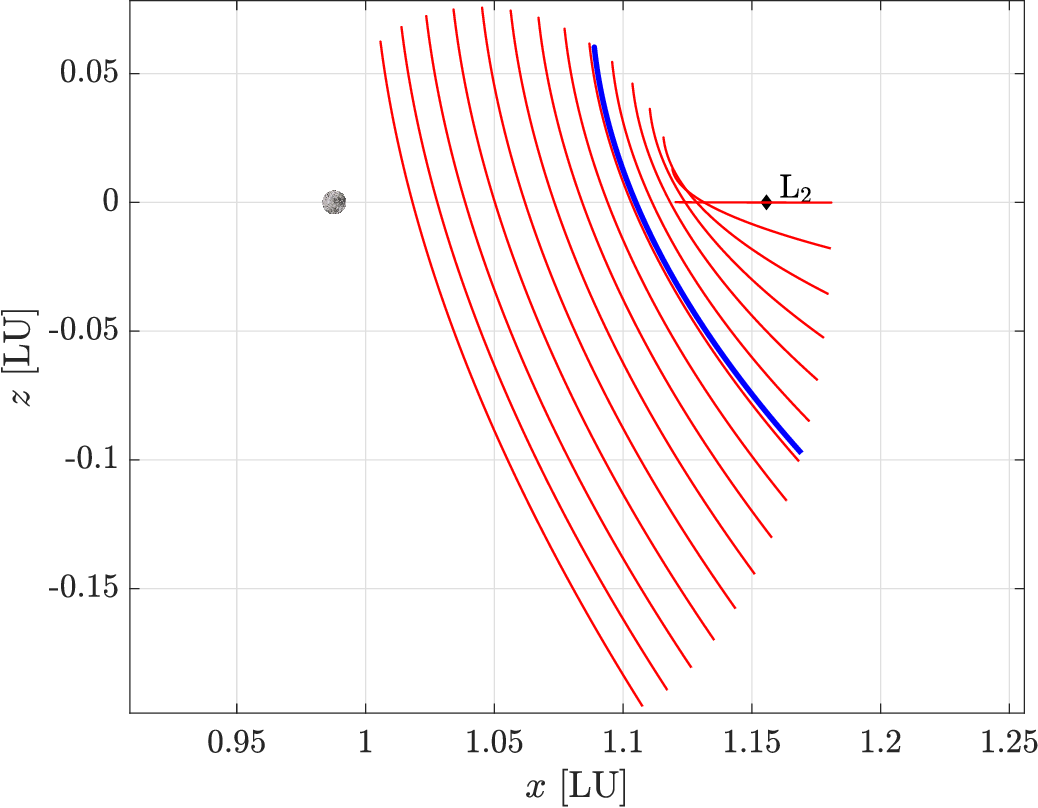}%
            \label{fig:phys-orbit-l2shalo_xz}
            }
        \subfloat[$yz$-plane.]{%
            \includegraphics[height=5.1em]{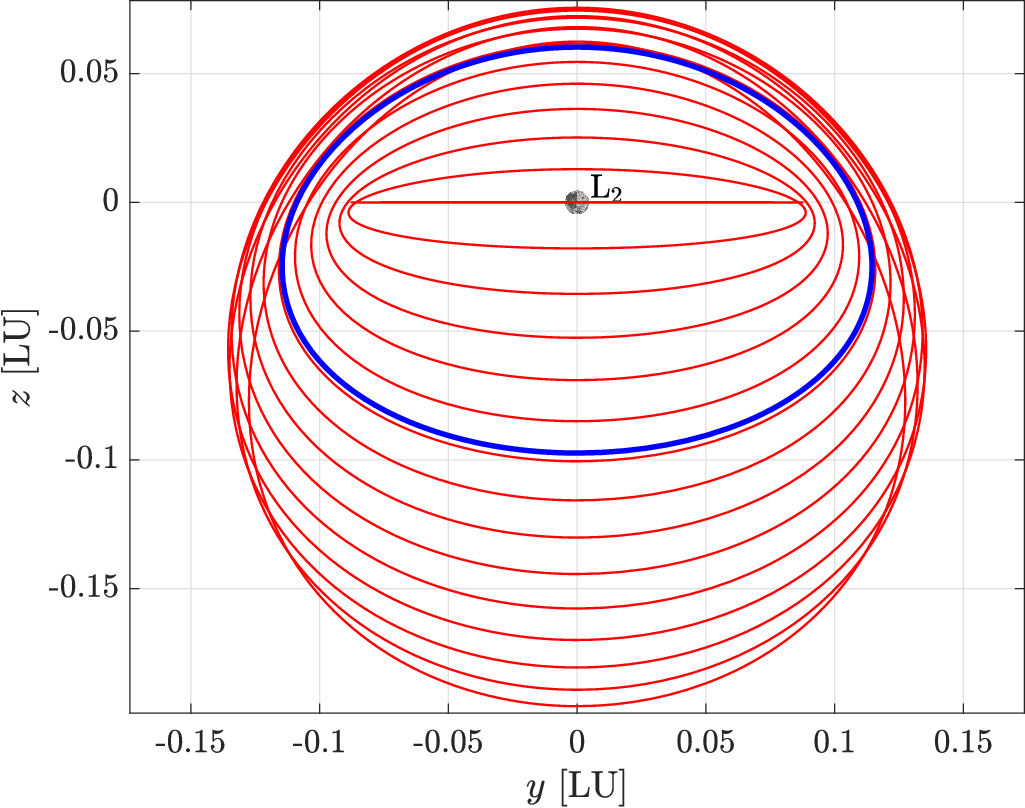}%
            \label{fig:phys-orbit-l2shalo_yz}
            }
        \hspace*{\fill}%
        \caption{Target $\tL_2$ southern halo orbit in the $\tL_2$ southern halo orbit family in the rotating reference frame with $\mfk{J}\geq{2.5}$.}
        \label{fig:phys-orbit-l2shalo}
    \end{minipage}
    \hspace*{\fill}%
    \begin{minipage}{0.48\textwidth}
        \centering
        \subfloat[3-dimensional.]{%
            \includegraphics[height=18.4em]{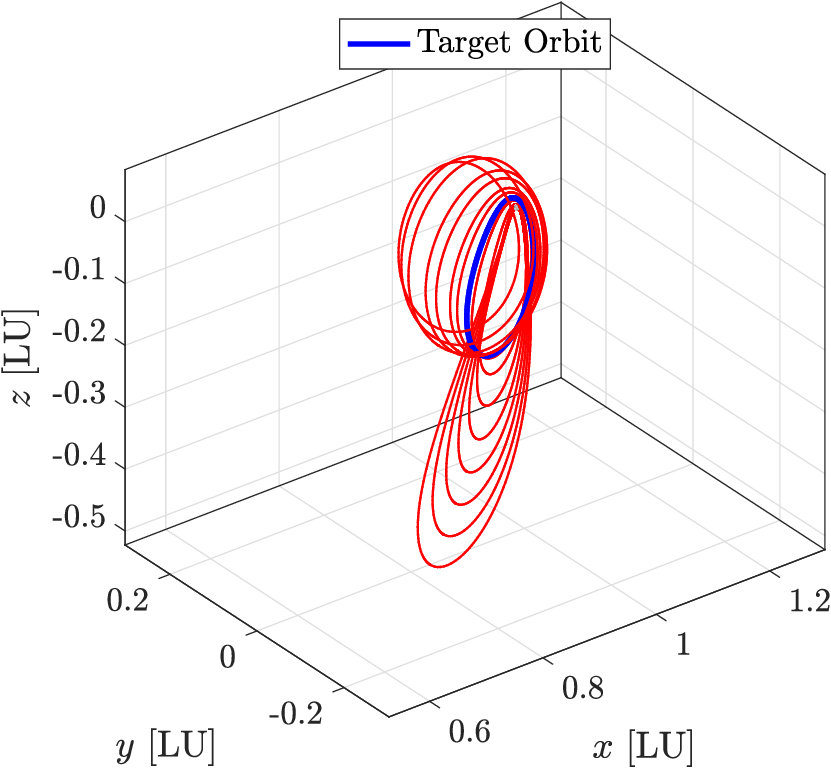}%
            \label{fig:phys-orbit-nrho_3d}
            }
        \\
        \hspace*{\fill}%
        \subfloat[$xy$-plane.]{%
            \includegraphics[height=5.1em]{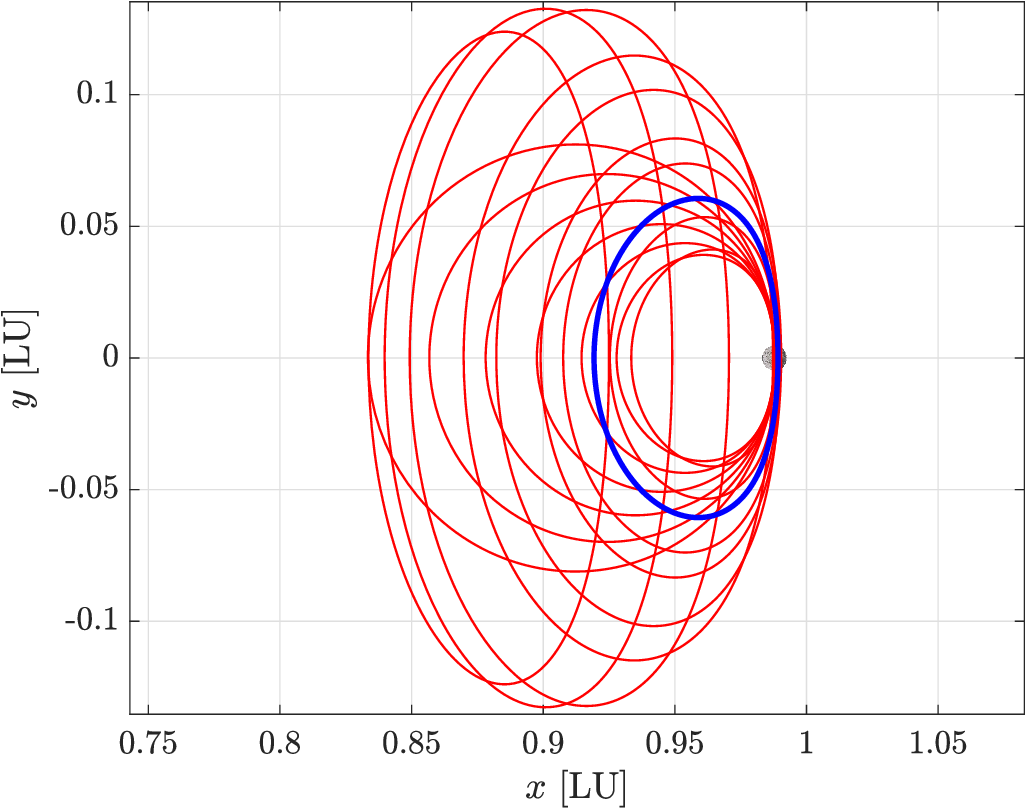}%
            \label{fig:phys-orbit-nrho_xy}
            }
        \subfloat[$xz$-plane.]{%
            \includegraphics[height=5.1em]{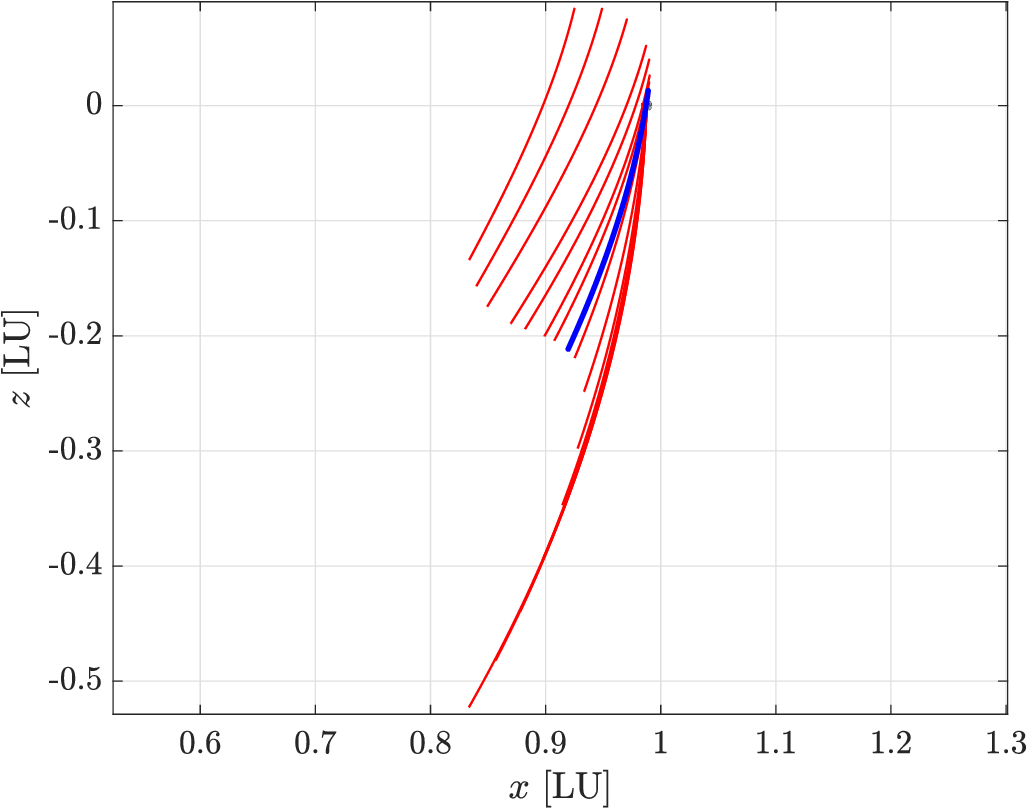}%
            \label{fig:phys-orbit-nrho_xz}
            }
        \subfloat[$yz$-plane.]{%
            \includegraphics[height=5.1em]{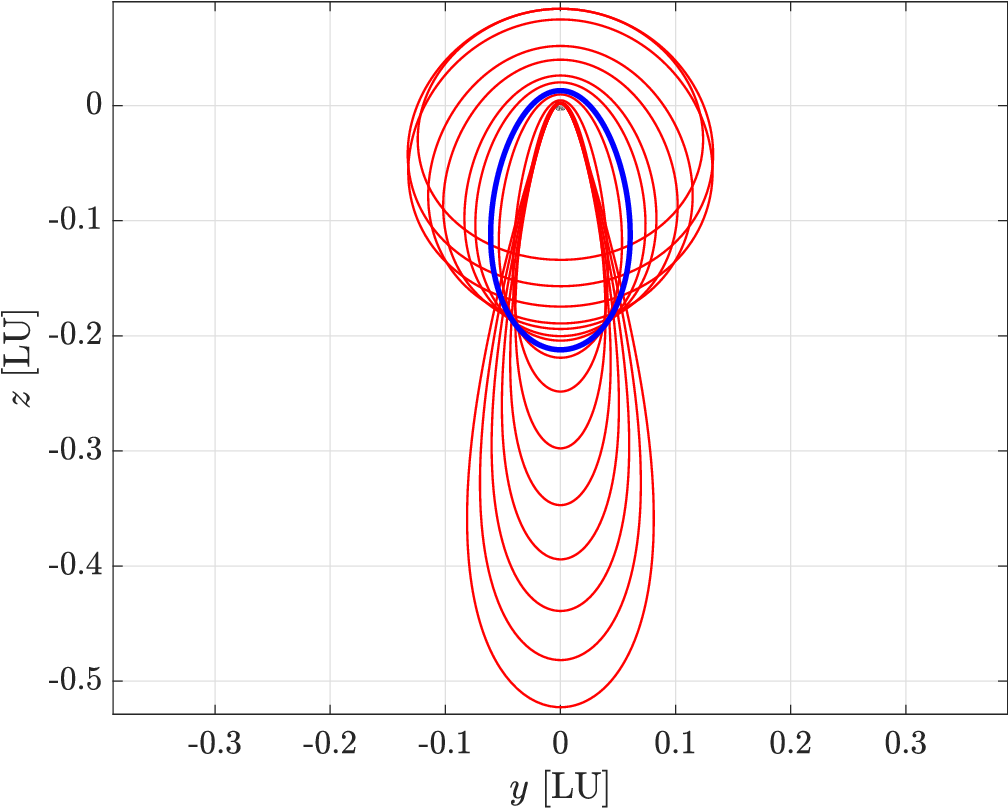}%
            \label{fig:phys-orbit-nrho_yz}
            }
        \hspace*{\fill}%
        \caption{Target near-rectilinear halo orbit in the $\tL_1$ southern halo orbit family in the rotating reference frame with $\mfk{J}\geq{2.7}$.}
        \label{fig:phys-orbit-nrho}
    \end{minipage}
\end{figure}

In this section, the physical orbit parameters utilized in the R3BP minimum-time, multi-mode optimal control problem formulation of Section~\ref{sec:prob} are presented.
In this investigation, the periodic orbits of interest are halo orbits in the Earth-Moon system.
The objective of this transfer is to demonstrate the trajectory generation procedure on a complex transfer in the Earth-Moon R3BP in order to solve the novel optimal control problem formulated in Section~\ref{sec:prob}.
The conditions defining the target periodic orbits are summarized in Table~\ref{tab:phys-orbit-prm}.
\begin{table}[!t]
    \centering
    \caption{Target CR3BP orbit parameters.}
    \begin{tabular*}{0.9\textwidth}{@{\extracolsep\fill}crrc@{}} \hline\hline
        \mrow{2}{\textbf{Parameter}} & \mcolc{1}{\textbf{Initial $\tL_2$ Southern}} & \mcolc{1}{\mrow{2}{\textbf{Terminal NRHO Value}}} & \mrow{2}{\textbf{Units}} \\
        & \mcolc{1}{\textbf{Halo Orbit Value}} & & \\ \hline
        $x$       &  1.1692032436399828$\times{10}^{+0~}$ &  9.1929792455210269$\times{10}^{-1~}$ & [LU] \\
        $y$       &  9.0895948914056935$\times{10}^{-30}$ & -7.8475765128918404$\times{10}^{-29}$ & [LU] \\
        $z$       & -9.7343078972773986$\times{10}^{-2~}$ & -2.1213093403317357$\times{10}^{-1~}$ & [LU] \\
        $\mad{x}$ & -1.2120988489044185$\times{10}^{-15}$ & -2.1943176228275019$\times{10}^{-13}$ & [LU/TU] \\
        $\mad{y}$ & -1.9424148423494397$\times{10}^{-1~}$ &  1.3779559700236524$\times{10}^{-1~}$ & [LU/TU] \\
        $\mad{z}$ &  1.1582044638613824$\times{10}^{-15}$ & -1.1031515999570366$\times{10}^{-12}$ & [LU/TU] \\
        $\mfk{J}$ &  3.1141257613953099$\times{10}^{+0~}$ &  3.0032754028672501$\times{10}^{+0~}$ & [LU$^2$/TU$^2$] \\
        $\mfk{p}$ &  3.3325377871055926$\times{10}^{+0~}$ &  1.8077163954358124$\times{10}^{+0~}$ & [TU] \\ \hline\hline
    \end{tabular*}
    \label{tab:phys-orbit-prm}
\end{table}
Specifically, Table~\ref{tab:phys-orbit-prm} presents the initial coast phase initial boundary condition defining the target 14.77-day $\tL_2$ southern halo orbit as well as the terminal coast phase terminal boundary condition defining the target 8.01-day NRHO (that is, within the $\tL_1$ southern halo orbit family).
The target orbits in their respective families are shown in Figures~\ref{fig:phys-orbit-l2shalo} and~\ref{fig:phys-orbit-nrho}.
Note that any CR3BP periodic orbit can be utilized in the problem formulation of Section~\ref{sec:prob}; however, the orbits shown in Table~\ref{subsec:phys-orbit} and Figures~\ref{fig:phys-orbit-l2shalo} and~\ref{fig:phys-orbit-nrho} are chosen to demonstrate a challenging, real-world, three-dimensional R3BP transfer within the lunar neighborhood.



\section{Numerical Approach}\label{sec:num}

The three-phase R3BP minimum-time, multi-mode periodic orbit transfer optimal control problem presented in Section~\ref{sec:prob} is solved using direct collocation.
Specifically, all results are obtained using the direct \textit{hp}-adaptive LGR collocation method~\cite{GargHuntington2010,GargRao2011b,PattersonRao2012} implemented in the \MATLAB{} general-purpose optimal control software \GPOPSII{}~\cite{PattersonRao2014} using the NLP solver SNOPT~\cite{GillSaunders2005}.
SNOPT is employed with a tolerance of $\epsilon_{\text{NLP}}=1\times{10}^{-8}$ and a maximum iteration count of 3000, where first derivatives are supplied via the \MATLAB{} automatic differentiation software ADiGator~\cite{WeinsteinRao2017}.
The LGR collocation method allows the optimal control problem to be constructed as a multiple-phase problem, where each phase contains separate dynamic, path, and boundary constraints.
Furthermore, the NLP arising from the discretization scheme is solved on a static mesh in each phase.
Because the size of the mesh (that is, the total number of discretization points in all phases) directly impacts the computational efficiency in which the NLP is solved, \textit{hp}-adaptive methods are employed to efficiently and accurately determine the placement of the discretization points.
The smooth mesh refinement method Ref.~\cite{HamanRao2025} is paired with the nonsmooth mesh refinement method of Ref.~\cite{MillerRao2021}, where a state error tolerance of $\epsilon=1\times{10}^{-6}$ and maximum number of mesh refinement iterations of $M_{\max}=25$ are employed.
The method of Ref.~\cite{HamanRao2025} is chosen here because it has been shown to efficiently determine an accurate static mesh that satisfies the desired state error tolerance, while simultaneously verifying the collocation solution with explicit propagation.
For use in the method of Ref.~\cite{HamanRao2025}, the following parameters are used: a minimum number of collocation points in an interval $N_{\min}=2$, maximum number of collocation points in an interval $N_{\max}=12$, and \MATLAB{} ODE solver \ode{113} employed with the options $\texttt{RelTol}=1\times{10}^{-10}$, $\texttt{AbsTol}=1\times{10}^{-10}$, and $\texttt{NormControl}=\texttt{on}$.
The method of Ref.~\cite{MillerRao2021} is chosen here because it accurately estimates locations of nonsmoothness in the optimal control.
For use in the method of Ref.~\cite{MillerRao2021}, the following parameters are used: a set of jump function approximation orders of $\mcl{R}=\{1,\ldots,8\}$, jump threshold of $\eta_{\tJ}=0.2$, and bound safety factor of $\eta_{\tB}=1.2$.
All periodic orbits are obtained using JPL's Three-Body Periodic Orbits toolbox~\cite{JPLSSD,VaqueroSenent2018}.
Finally, all computations are performed on a 12-core Apple M3 Pro MacBook Pro running macOS Sonoma version 14.6.1 with 36 GB LPDDR5 of unified memory using \MATLAB{} version R2024b (24.2.0.2740171 Update 1).

The solution generation procedure follows similar steps to that in a continuation method, where the solution to a neighboring (and often less complex) problem is solved and used as the initial guess for the problem of interest.
First, a propellant-unconstrained, minimum-time solution is obtained using only mode 1 (that is, the high-thrust mode), where no upper bound is placed on the allowable amount of mode 1 propellant consumed.
The propellant-unconstrained, minimum-time transfer between the two specified periodic orbits using only mode 1 is obtained using an initial guess obtained via \textit{trajectory stacking}, which is straightforward, requires little intuition the problem, and increases computational efficiency relative to a straight-line initial guess.
The trajectory stacking method proceeds as follows~\cite{HamanRao2024a,HamanRao2024b}.
The initial orbit is propagated forward in time from the specified initial condition for a specified number of orbital periods (that is, stacking the initial orbit on top of itself).
Then, the terminal orbit is propagated backward in time from the specified terminal condition for a specified number of orbital periods (that is, stacking the terminal orbit on top of itself).
The two previously generated trajectories are then combined to form the initial guess trajectory, which represents a feasible, ballistic solution with the exception of the state discontinuity at the epoch where the two trajectories are patched together.
Upon analyzing the amount of mode 1 propellant required for the propellant-unconstrained, only mode 1 minimum-time solution, a constraint is placed on the allowable amount of mode 1 propellant consumed.
Upon making no assumptions regarding the optimal control structure, the problem is then re-solved using only mode 1, only mode 2, or multi-mode propulsion capabilities with a specified propellant constraint value, where the unconstrained solution (that is, using the associated true anomaly, state, control, static parameter, and mesh solutions) is used as the initial guess.
Once the optimal control structure is determined for the transfer phase, the transfer phase is partitioned into multiple domains in order to optimize the throttle switch times, which adopts a similar approach for optimizing bang-bang control structures as developed in Ref.~\cite{PagerRao2022}.
Finally, the propellant constraint value is iteratively decreased, where the solution obtained using the previous propellant constraint value is used as the initial guess.


\section{Results and Discussion}\label{sec:res}

Using the physical parameters and numerical approach of Sections~\ref{sec:phys} and~\ref{sec:num}, respectively, this section presents the results acquired by solving the three-phase R3BP minimum-time, multi-mode periodic orbit transfer optimal control problem of Section~\ref{sec:prob}.
Because the libration point orbits of interest are periodic solutions, a control effort is required for the spacecraft to depart the initial orbit and arrive on the terminal orbit; thus, it is expected that there are at least two thrusting regions: (1) departing the initial orbit and (2) arriving on the terminal orbit.
Furthermore, while the invariant manifolds associated with the periodic orbits can be exploited for energy-efficient and minimum-fuel transfers, it is expected that the minimum-time transfers utilize neither the unstable invariant manifold associated with the initial periodic orbit nor the stable invariant manifold associated with the terminal periodic orbit.

The remainder of this section is organized as follows.
First, Section~\ref{subsec:res-base} discusses the baseline results that determine the optimal control structure using (1) only mode 1 and (2) multi-mode propulsion.
Then, Section~\ref{subsec:res-m1} presents the results obtained using only mode 1 propulsion, while Section~\ref{subsec:res-mm} presents the results obtained using multi-mode propulsion.
Finally, Section~\ref{subsec:res-sum} summarizes the results.

\subsection{Baseline Results}\label{subsec:res-base}

\begin{figure}[!b]
    \centering
    \subfloat[3-dimensional.]{%
        \includegraphics[height=16.3em]{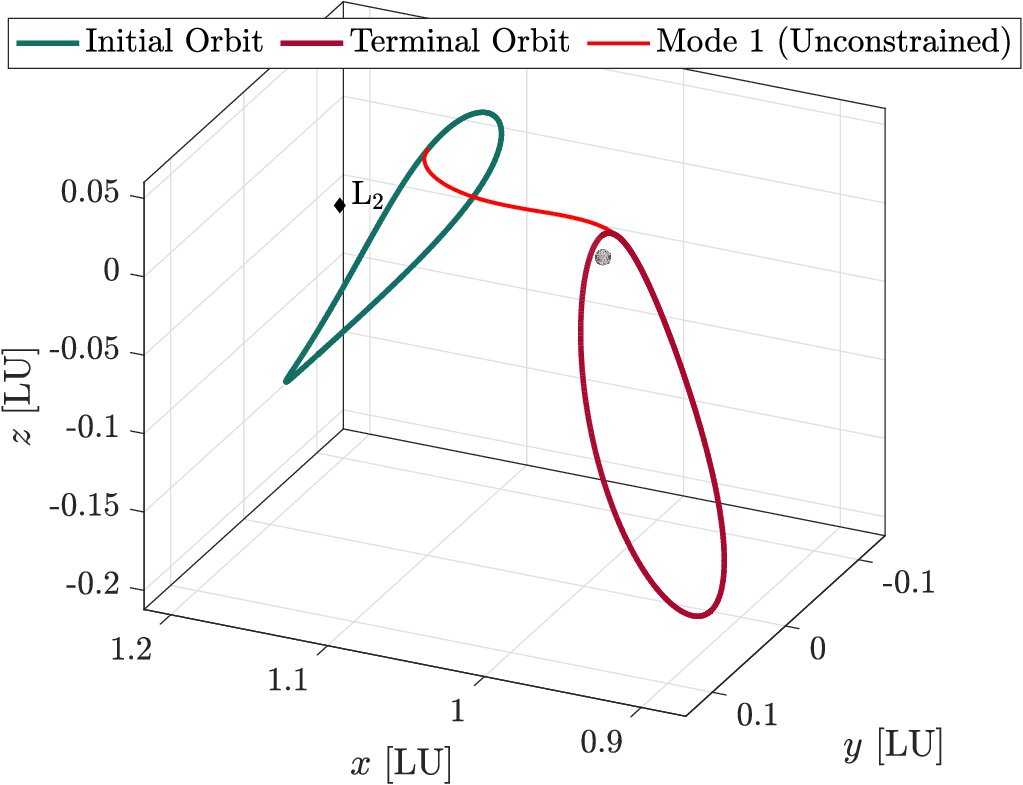}%
        \label{fig:res-base-m1f_3d}
        }
    \\
    \hspace*{\fill}%
    \subfloat[$xy$-plane.]{%
        \includegraphics[height=5.3em]{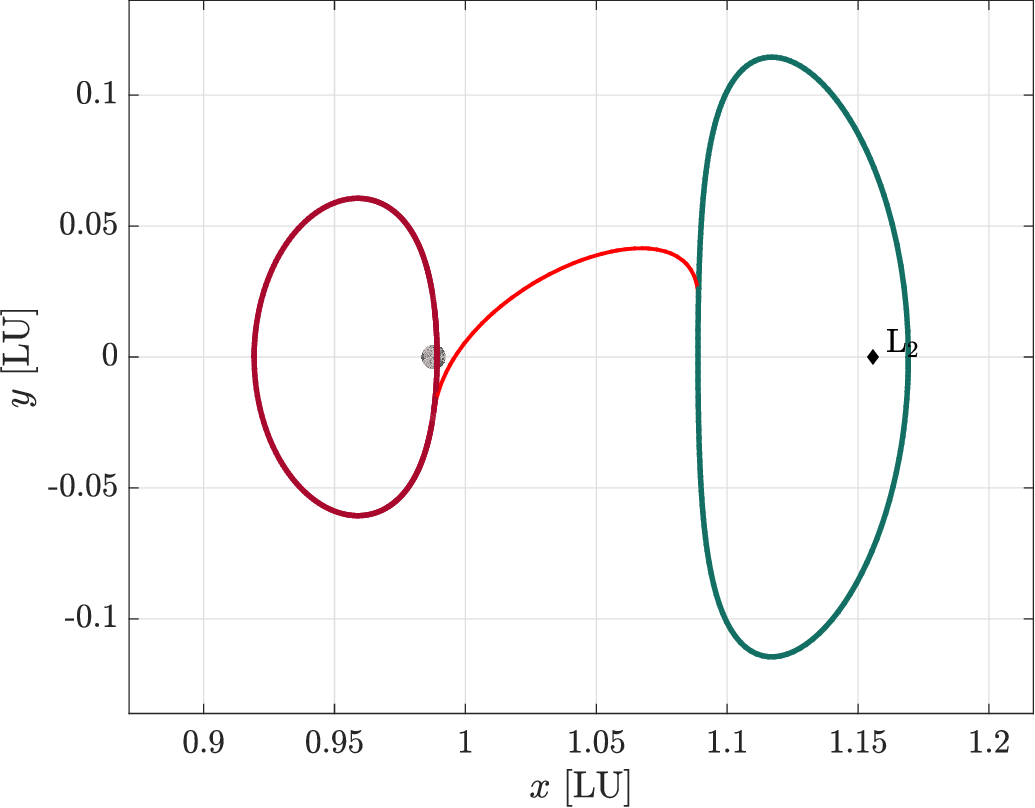}%
        \label{fig:res-base-m1f_xy}
        }
    \subfloat[$xz$-plane.]{%
        \includegraphics[height=5.3em]{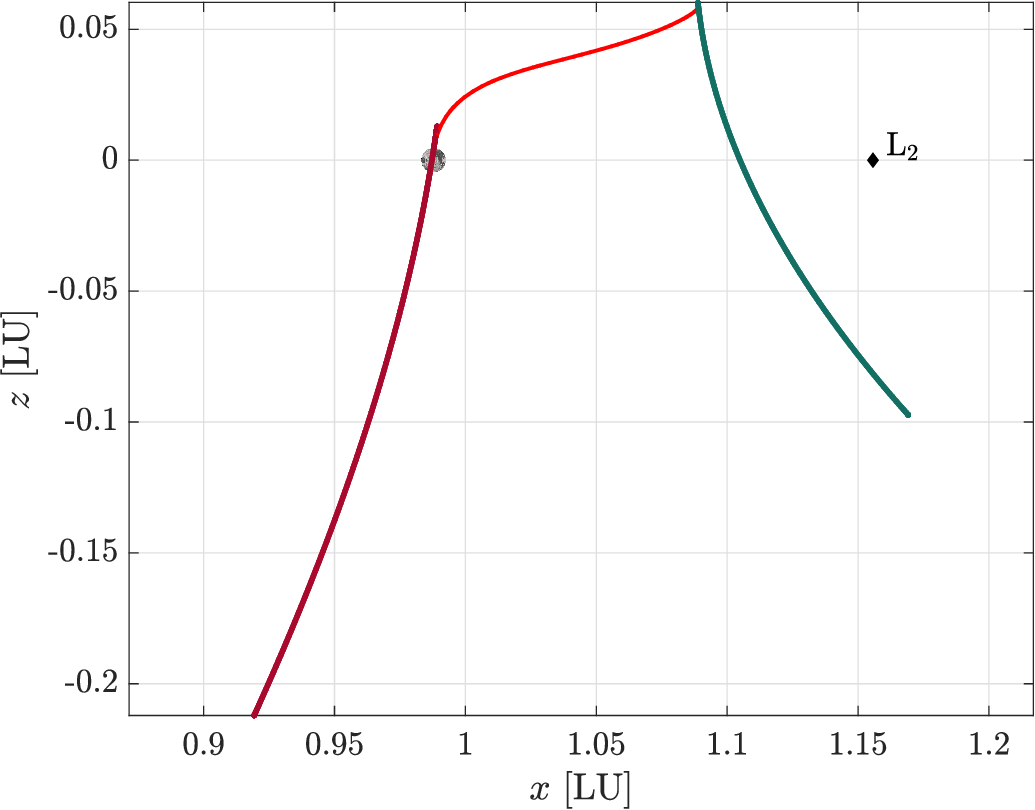}%
        \label{fig:res-base-m1f_xz}
        }
    \subfloat[$yz$-plane.]{%
        \includegraphics[height=5.3em]{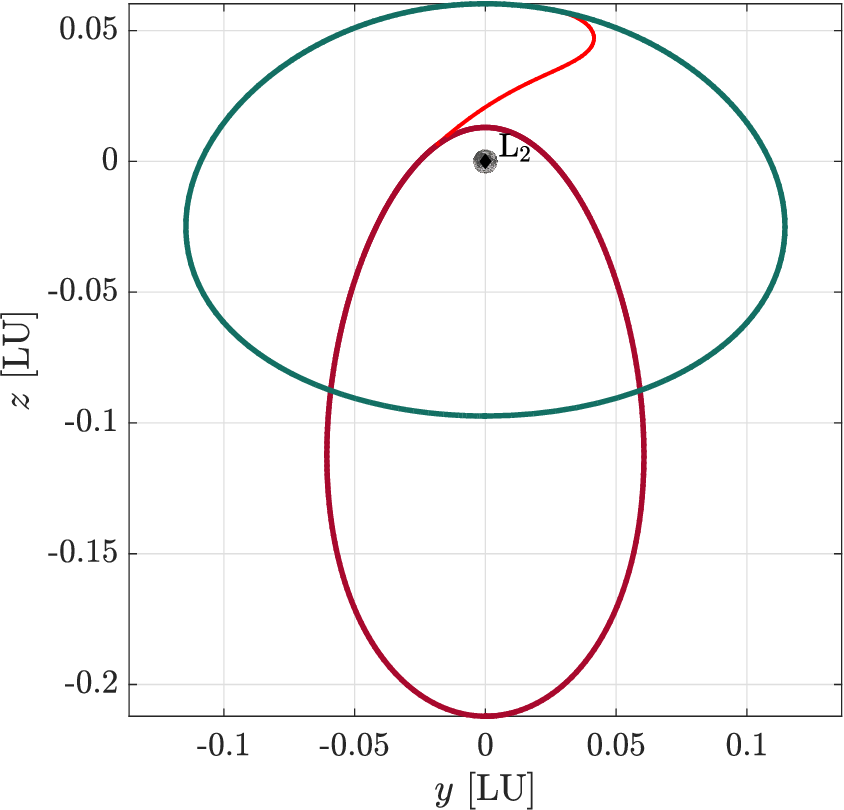}%
        \label{fig:res-base-m1f_yz}
        }
    \hspace*{\fill}%
    \caption{Propellant-unconstrained, minimum-time transfer using only mode 1 propulsion.}
    \label{fig:res-base-m1f}
\end{figure}

\begin{figure}[!t]
    \centering
    \hspace*{\fill}%
    \includegraphics[height=15.0em]{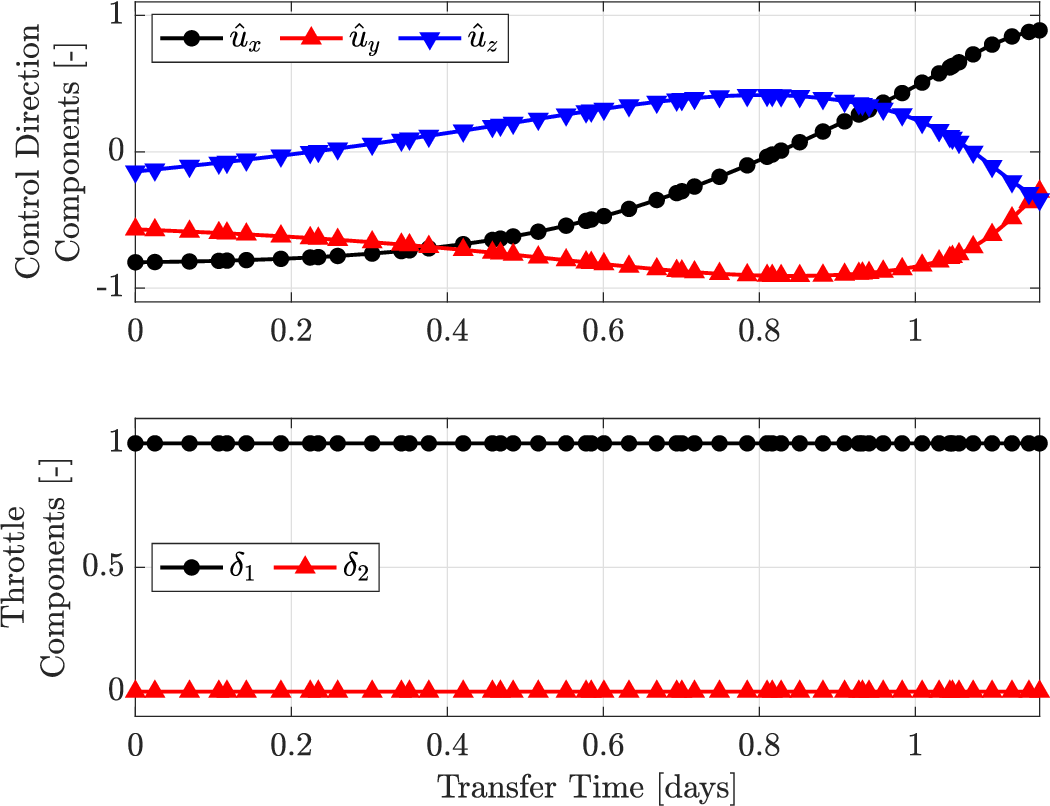}%
    \hspace*{\fill}%
    \caption{Optimal control direction and throttle component history corresponding to Figure~\ref{fig:res-base-m1f}.}
    \label{fig:res-base-m1f_u}
\end{figure}

In this section, the three-phase minimum-time problem is solved using only mode 1 propulsion to provide a baseline trajectory and control structure.
For the propellant-unconstrained problem, an initial guess is constructed via trajectory stacking, where $\delta_1\in[0,1]$ and $\delta_2=0$.
Using this initial guess, the propellant-unconstrained, minimum-time transfer is shown in Figure~\ref{fig:res-base-m1f}, and the corresponding optimal control direction and throttle component profiles are shown in Figure~\ref{fig:res-base-m1f_u}.%
For this solution, an optimal objective value of $\mcl{J}^{*}=0.285471$ (that is, 1.163 [days]) is obtained, where mode 1 consumes 40.973 [kg] of fuel.
From the $\tL_2$ southern halo orbit initial condition defined in Table~\ref{tab:phys-orbit-prm}, the spacecraft coasts along the initial orbit for 7.751 [days] (that is, approximately 52.471\% of the orbital period) before departing the initial orbit.
Upon arriving on the NRHO, the spacecraft then coasts along the terminal orbit for 3.921 [days] (that is, approximately 48.928\% of the orbital period) to the terminal condition defined in Table~\ref{tab:phys-orbit-prm}.
As shown in Figure~\ref{fig:res-base-m1f_u}, the control direction component profiles are generally smooth, and the optimal mode 1 throttle is at its maximum (that is, $\delta_1=1$) for the entire transfer phase, which is typical for propellant-unconstrained, minimum-time solutions.

\begin{figure}[!t]
    \centering
    \begin{minipage}{0.48\textwidth}
        \centering
        \hspace*{\fill}%
        \includegraphics[height=15.0em]{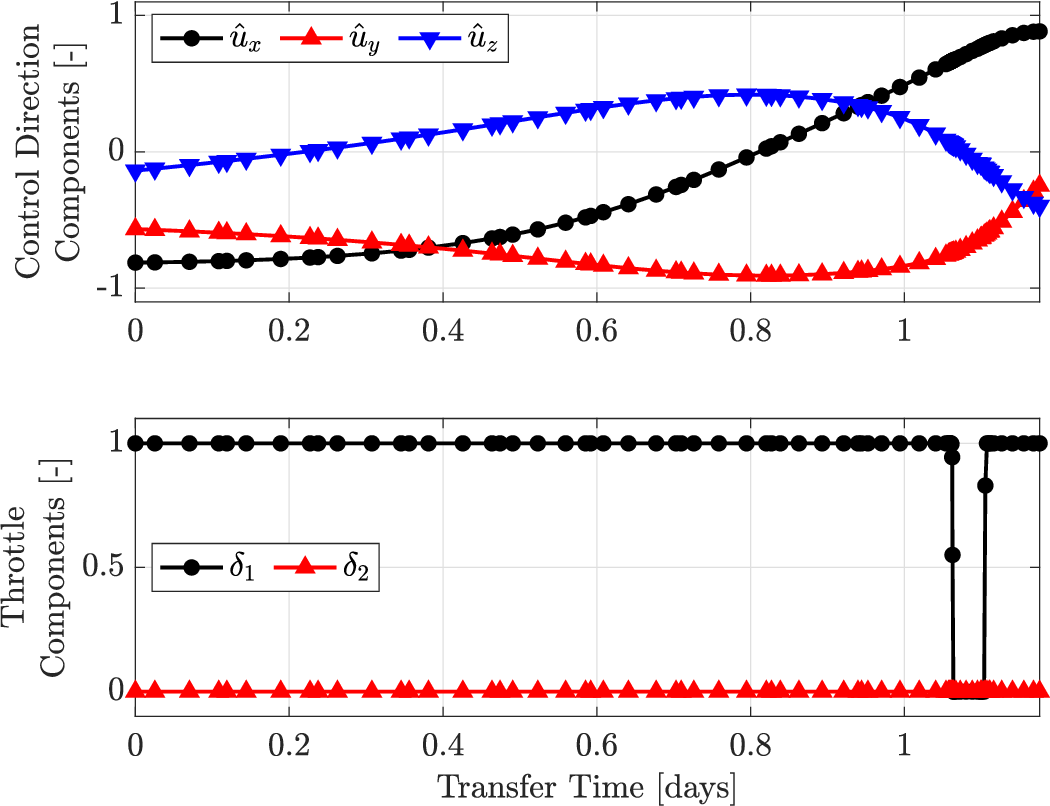}%
        \hspace*{\fill}%
        \caption{Baseline optimal control direction and throttle component history for the 40 [kg] propellant-constrained, minimum-time transfer using only mode 1 propulsion.}
        \label{fig:res-base-m1b_u}
    \end{minipage}
    \hspace*{\fill}%
    \begin{minipage}{0.48\textwidth}
        \hspace*{\fill}%
        \includegraphics[height=15.0em]{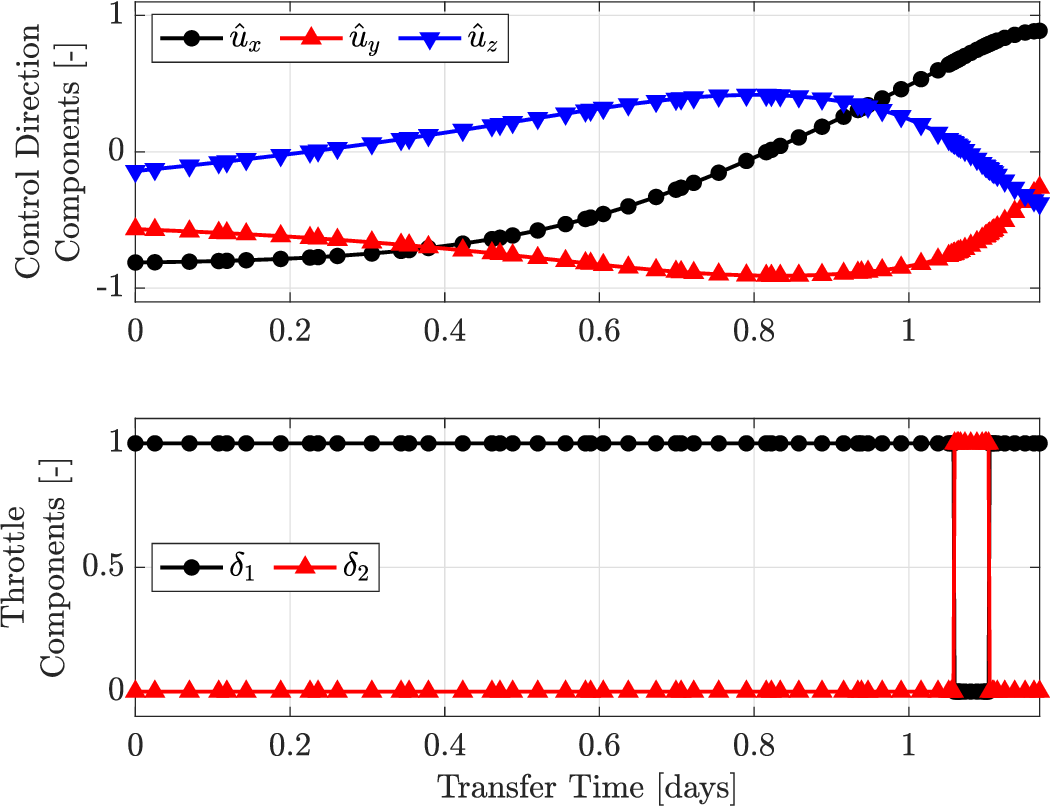}%
        \hspace*{\fill}%
        \caption{Baseline optimal control direction and throttle component history for the 40 [kg] propellant-constrained, minimum-time transfer using multi-mode propulsion ($T_{\max,2}=0.5$ [N]).}
        \label{fig:res-base-mmb_ts1_u}
    \end{minipage}
\end{figure}

Now that the propellant-unconstrained solution provides a baseline value for the amount of propellant consumed by mode 1, propellant constraint values less than 40.973 [kg] are enforced.
As previously stated, it is expected that mode 1 turns off for a finite duration due to the propellant constraint, which results in switches in the mode 1 control structure.
Then, using the propellant-unconstrained, only mode 1 solution as the initial guess, the propellant-constrained, minimum-time problem is solved for two cases: (1) only mode 1 and (2) multi-mode.
For a propellant constraint value of 40 [kg], the optimal control direction and throttle component histories are shown in Figures~\ref{fig:res-base-m1b_u} and~\ref{fig:res-base-mmb_ts1_u} using only mode 1 and multi-mode with $T_{\max,2}=0.5$ [N], respectively.
Upon making no assumptions regarding the optimal control structure, the mesh refinement methods of Refs.~\cite{HamanRao2025} and~\cite{MillerRao2021} determine the locations of the switches in the control structure by increasing the mesh density in the respective regions.
As shown in Figure~\ref{fig:res-base-m1b_u}, the optimal control structure for the only mode 1 solution is determined to be on-off-on.
On the other hand as shown in Figure~\ref{fig:res-base-mmb_ts1_u}, the optimal control structure for the multi-mode solution is determined to be on-off-on for mode 1 and off-on-off for mode 2.
A similar structure to Figure~\ref{fig:res-base-mmb_ts1_u} is observed for the multi-mode case with $T_{\max,2}=0.25$ [N].
While the discontinuity detection method of Ref.~\cite{MillerRao2021} seems sufficient enough to determine the control structure on its own, lower propellant constraint values often cause (1) more switches in the control structure to be identified and (2) the propellant constraint value to not be reached (that is, a suboptimal solution is obtained).
In fact, the solution corresponding to Figure~\ref{fig:res-base-mmb_ts1_u} only utilizes 39.727 [kg] of the full 40 [kg] of fuel allotted to mode 1, which suggests a lower transfer time is possible.
To help avoid such complications, the transfer phase for the propellant-constrained, minimum-time problem is partitioned into domains based on the control structures observed in Figures~\ref{fig:res-base-m1b_u} and~\ref{fig:res-base-mmb_ts1_u} in order to optimize the location of the control switch times.
The procedure and results obtained using only mode 1 and multi-mode propulsion are discussed further in Sections~\ref{subsec:res-m1} and~\ref{subsec:res-mm}, respectively.


\subsection{Only Mode 1 Results}\label{subsec:res-m1}

\begin{figure}[!t]
    \centering
    \subfloat[3-dimensional.]{%
        \includegraphics[height=16.3em]{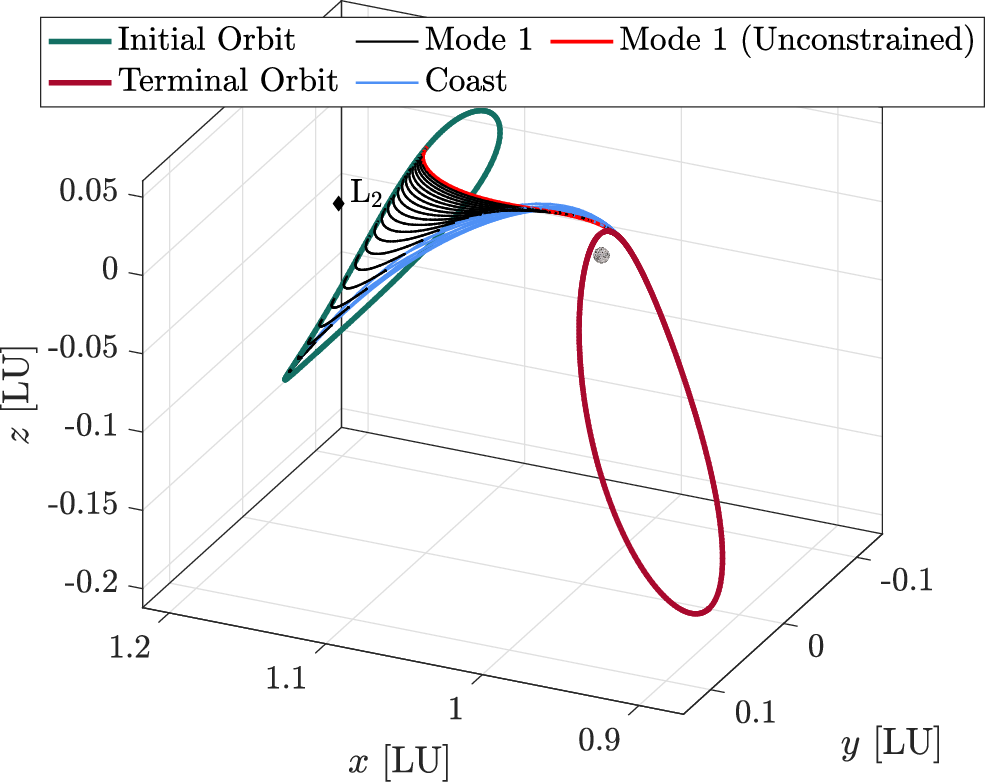}%
        \label{fig:res-m1-3d}
        }
    \\
    \hspace*{\fill}%
    \subfloat[$xy$-plane.]{%
        \includegraphics[height=5.3em]{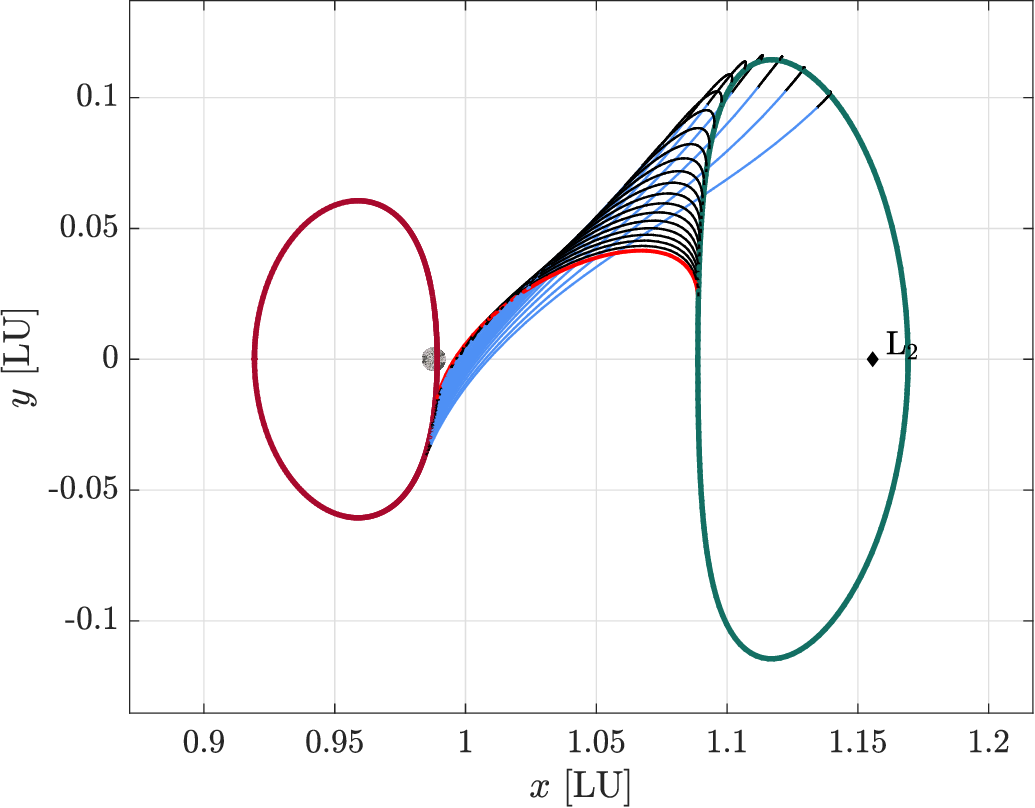}%
        \label{fig:res-m1-xy}
        }
    \subfloat[$xz$-plane.]{%
        \includegraphics[height=5.3em]{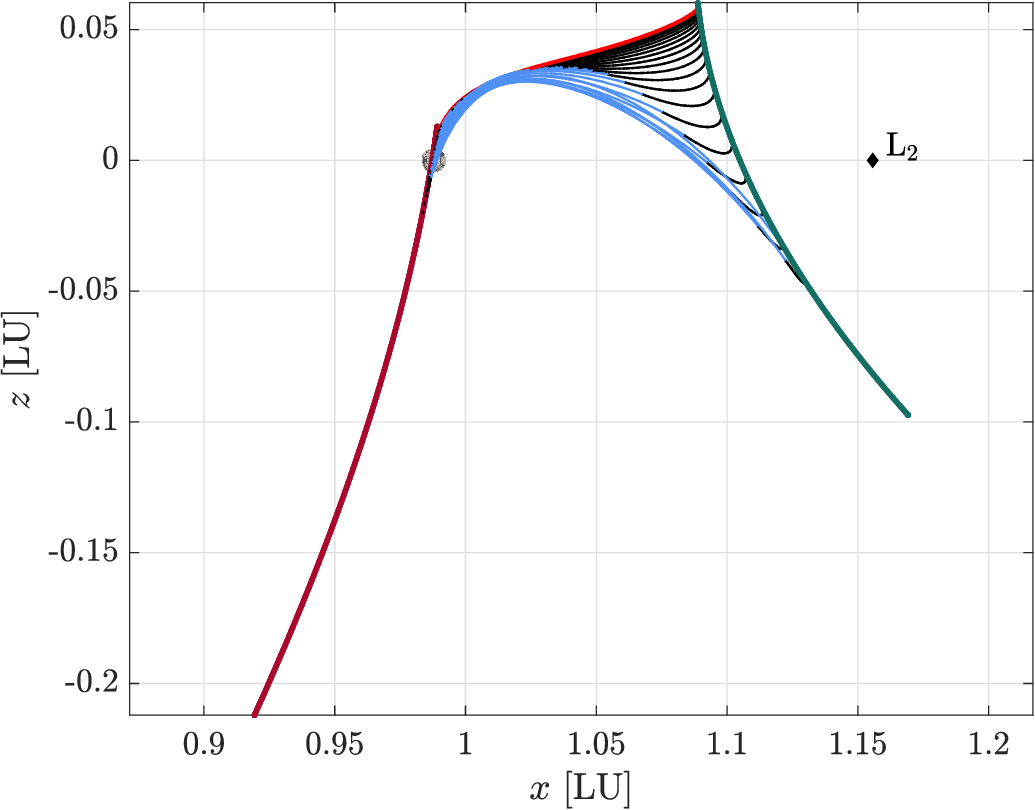}%
        \label{fig:res-m1-xz}
        }
    \subfloat[$yz$-plane.]{%
        \includegraphics[height=5.3em]{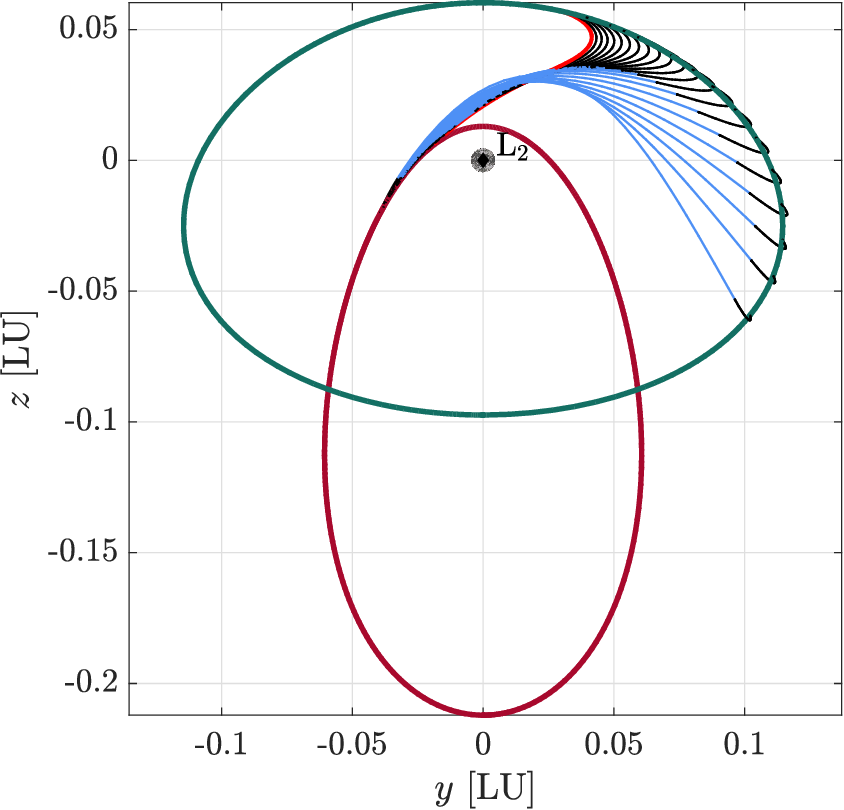}%
        \label{fig:res-m1-yz}
        }
    \hspace*{\fill}%
    \caption{Propellant-constrained, minimum-time transfers using only mode 1 propulsion.}
    \label{fig:res-m1-m1}
\end{figure}

The three-phase propellant-constrained, minimum-time problem is solved for propellant constraint values of 40 [kg] and less using only mode 1 propulsion with the optimal control structure determined in Section~\ref{subsec:res-base}.
As a result, the transfer phase is partitioned into 3 domains to optimize the throttle switch times.
The throttle structure for each domain is as follows: (1) mode 1: on, mode 2: off, (2) mode 1: off, mode 2: off, and (3) mode 1: on, mode 2: off.
Note that an active mode is equivalent to full throttle.
Using the numerical approach of Section~\ref{sec:num}, the trajectory envelope is shown in Figure~\ref{fig:res-m1-m1}, where the mode 1 thrusting regions are highlighted in black and the coasting regions are highlighted in pale blue.
For reference, the propellant-unconstrained mode 1 only solution shown in Figure~\ref{fig:res-base-m1f} is highlighted in red in Figure~\ref{fig:res-m1-m1}.

For the results shown in Figure~\ref{fig:res-m1-m1}, propellant constraint values from 40 [kg] to 20 [kg] are analyzed.
The 40 [kg] propellant-constrained trajectory is closest to the propellant-unconstrained solution using only mode 1, and for this solution, an optimal objective value of $\mcl{J}^{*}=0.289159$ (that is, 1.178 [days]) is obtained.
The optimal mode 1 control profile consists of the following: (1) a 1.063 [days] thrusting arc consuming 37.451 [kg] of fuel, (2) a 0.043 [days] coast arc, and (3) a 0.072 [days] thrusting arc consuming 2.549 [kg] of fuel.
As the propellant constraint value is decreased, the optimal trajectory deviates further and further from this solution, as shown in Figure~\ref{fig:res-m1-m1}.
For the 20 [kg] propellant-constrained trajectory, an optimal objective value of $\mcl{J}^{*}=0.632199$ (that is, 2.574 [days]) is obtained.
For this solution, the optimal mode 1 control profile consists of the following: (1) a 0.445 [days] thrusting arc consuming 15.683 [kg] of fuel, (2) a 2.007 [days] coast arc, and (3) a 0.122 [days] thrusting arc consuming 4.317 [kg] of fuel.

\begin{figure}[!t]
    \centering
    \hspace*{\fill}%
    \includegraphics[height=16.0em]{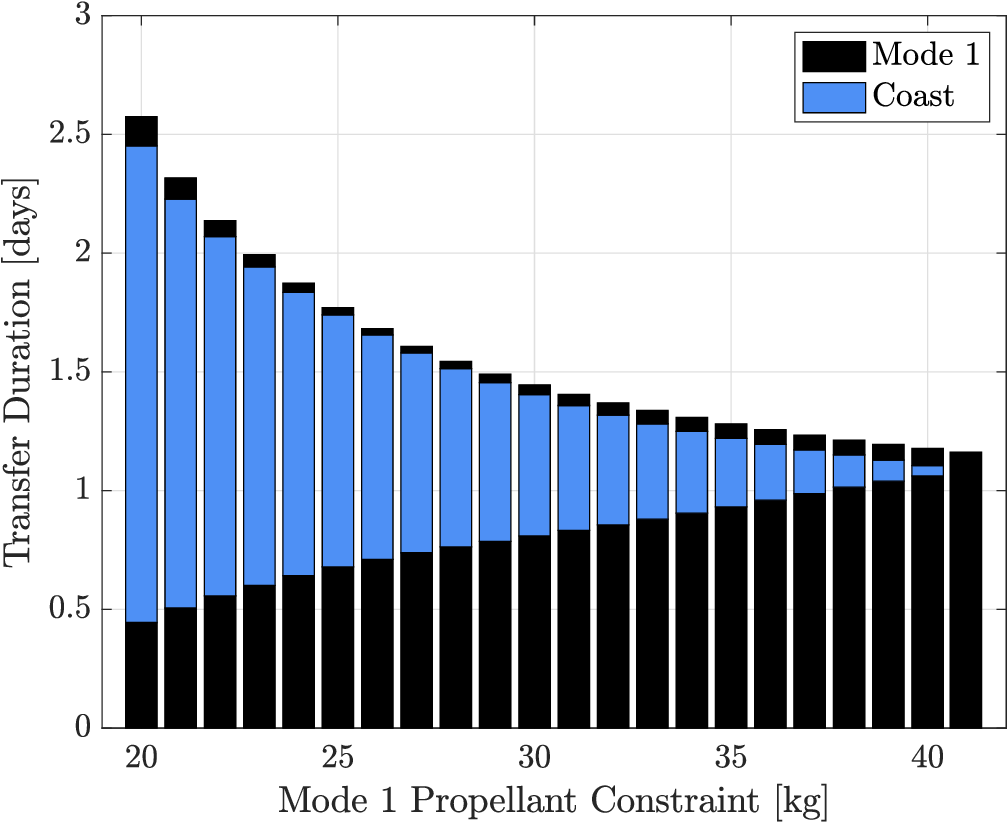}%
    \hspace*{\fill}%
    \caption{Transfer phase thrusting history for transfers shown in Figure~\ref{fig:res-m1-m1}.}
    \label{fig:res-m1-tpdursplit}
\end{figure}

A summary of the optimal mode 1 control profile for the transfer phase corresponding to Figure~\ref{fig:res-m1-m1} is shown in Figure~\ref{fig:res-m1-tpdursplit}.
As mentioned in Section~\ref{subsec:res-base}, the propellant-unconstrained, minimum-time transfer using only mode 1 propulsion consumes 40.973 [kg] of fuel, which is shown on the far right in Figure~\ref{fig:res-m1-tpdursplit} with no coast arc.
As the propellant constraint value is decreased from 40 [kg] to 20 [kg], the duration of the coast arc increases from 3.619\% to 77.958\% of the transfer phase duration, respectively.
While the duration of the thrusting arc arriving on the terminal orbit does not exhibit a clear trend, the duration of the thrusting arc departing from the initial orbit decreases from 90.238\% to 17.284\% of the transfer phase duration, respectively.
For all propellant constraint values ranging from 40 [kg] to 20 [kg], the optimal mode 1 throttle structure of on-off-on (that is, with two control switches) is observed, where no convergence issues are encountered.
For a propellant constraint value of 19 [kg], convergence is not achieved; thus, it is expected that the minimum-fuel solution using only mode 1 propulsion utilizes 19-20 [kg] of fuel.


\subsection{Multi-Mode Results}\label{subsec:res-mm}

\begin{figure}[!t]
    \centering
    \begin{minipage}{0.48\textwidth}
        \centering
        \subfloat[3-dimensional.]{%
            \includegraphics[height=16.3em]{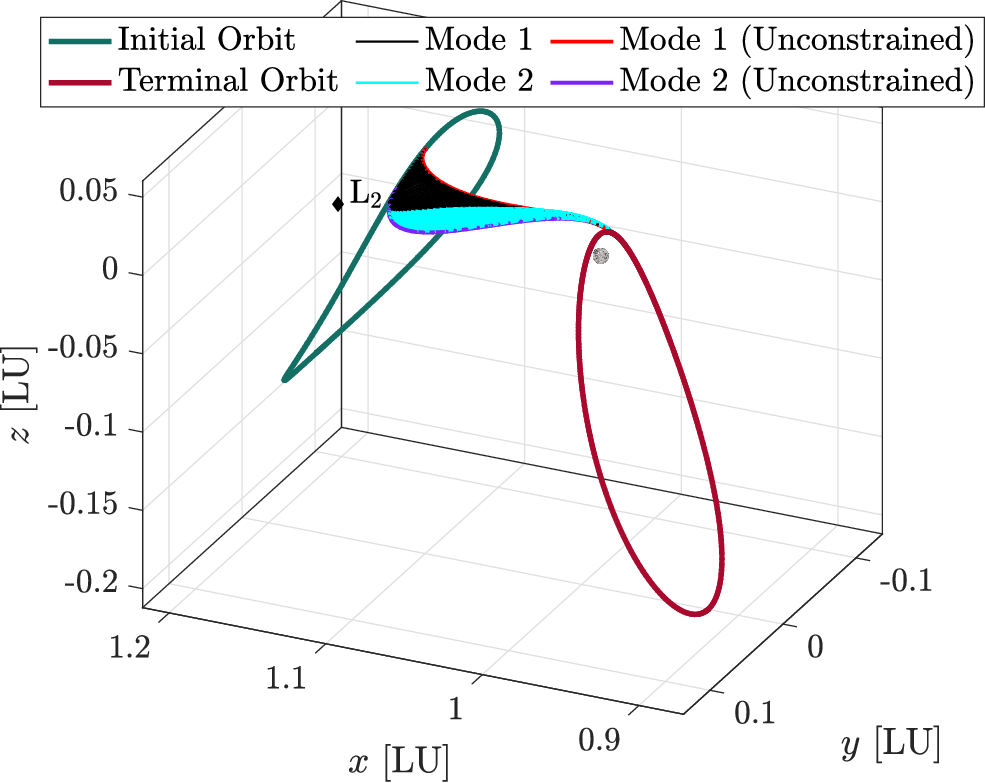}%
            \label{fig:res-mm-ts1_3d}
            }
        \\
        \hspace*{\fill}%
        \subfloat[$xy$-plane.]{%
            \includegraphics[height=5.3em]{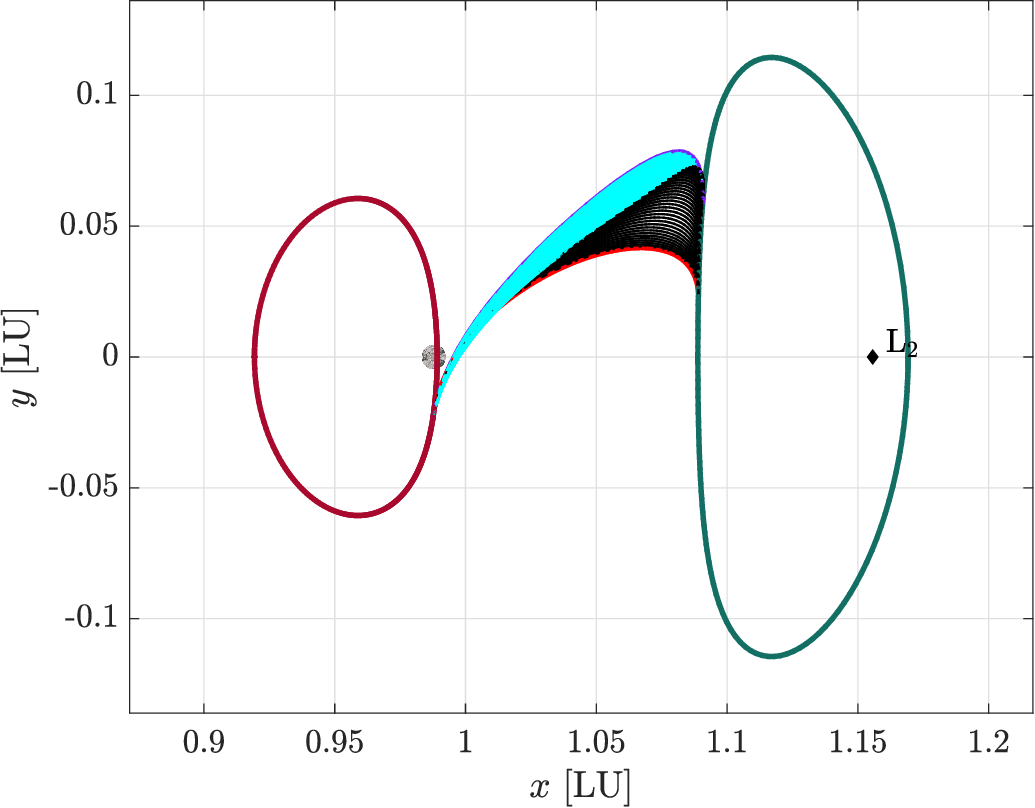}%
            \label{fig:res-mm-ts1_xy}
            }
        \subfloat[$xz$-plane.]{%
            \includegraphics[height=5.3em]{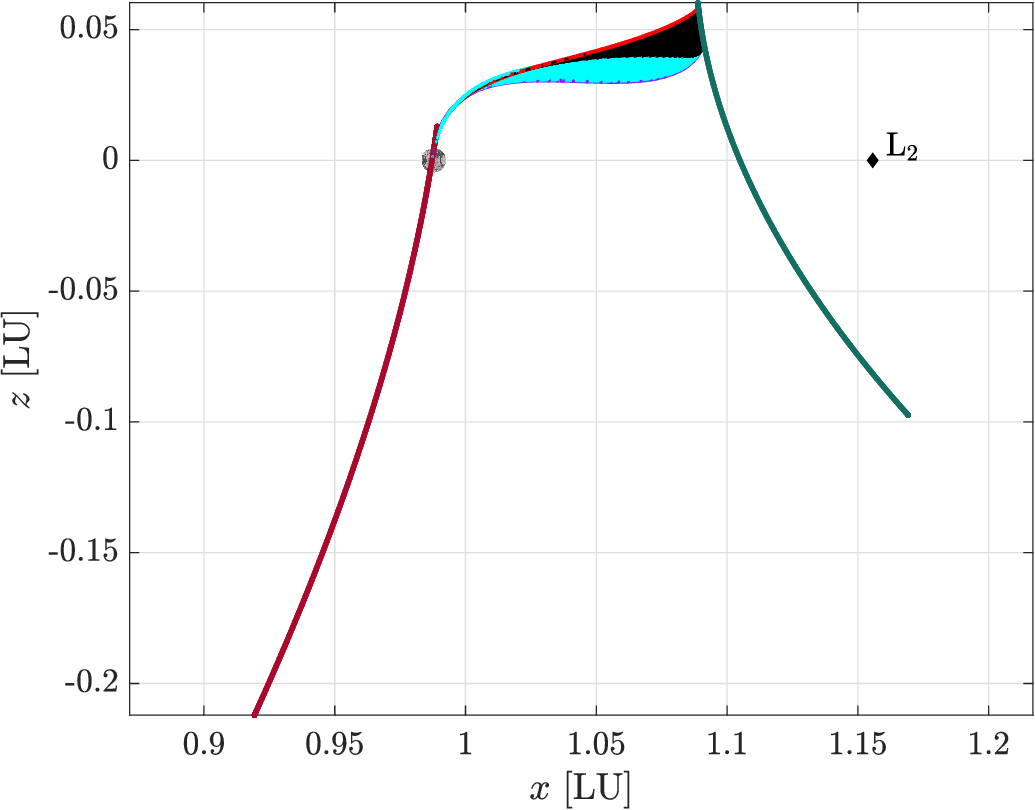}%
            \label{fig:res-mm-ts1_xz}
            }
        \subfloat[$yz$-plane.]{%
            \includegraphics[height=5.3em]{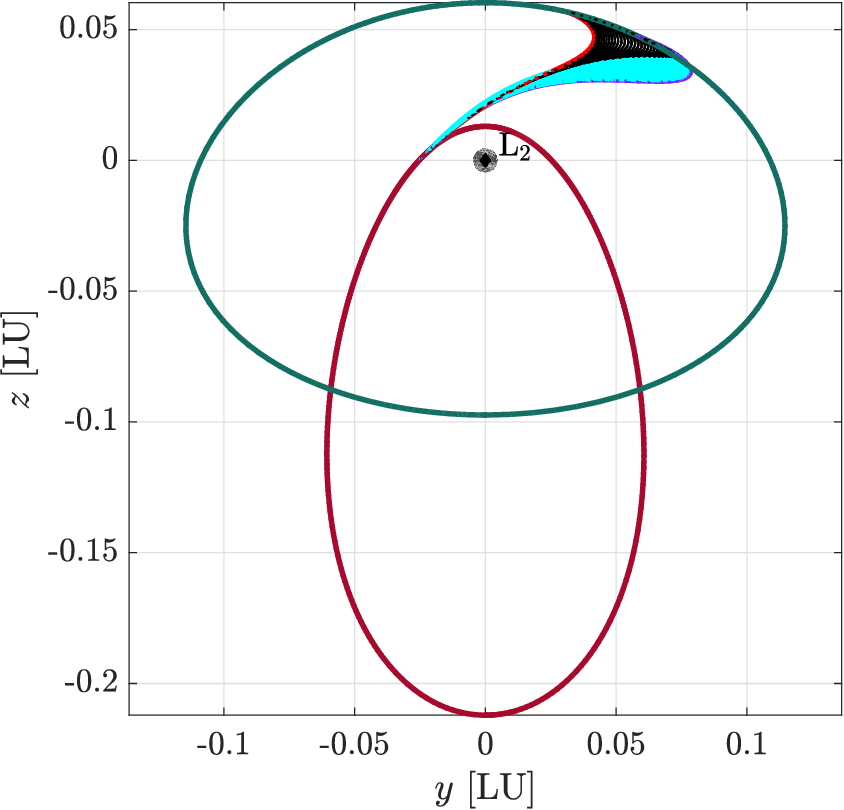}%
            \label{fig:res-mm-ts1_yz}
            }
        \hspace*{\fill}%
        \caption{Propellant-constrained, minimum-time transfers using multi-mode propulsion ($T_{\max,2}=0.50$ [N]).}
        \label{fig:res-mm-ts1}
    \end{minipage}
    \hspace*{\fill}%
    \begin{minipage}{0.48\textwidth}
        \centering
        \subfloat[3-dimensional.]{%
            \includegraphics[height=16.3em]{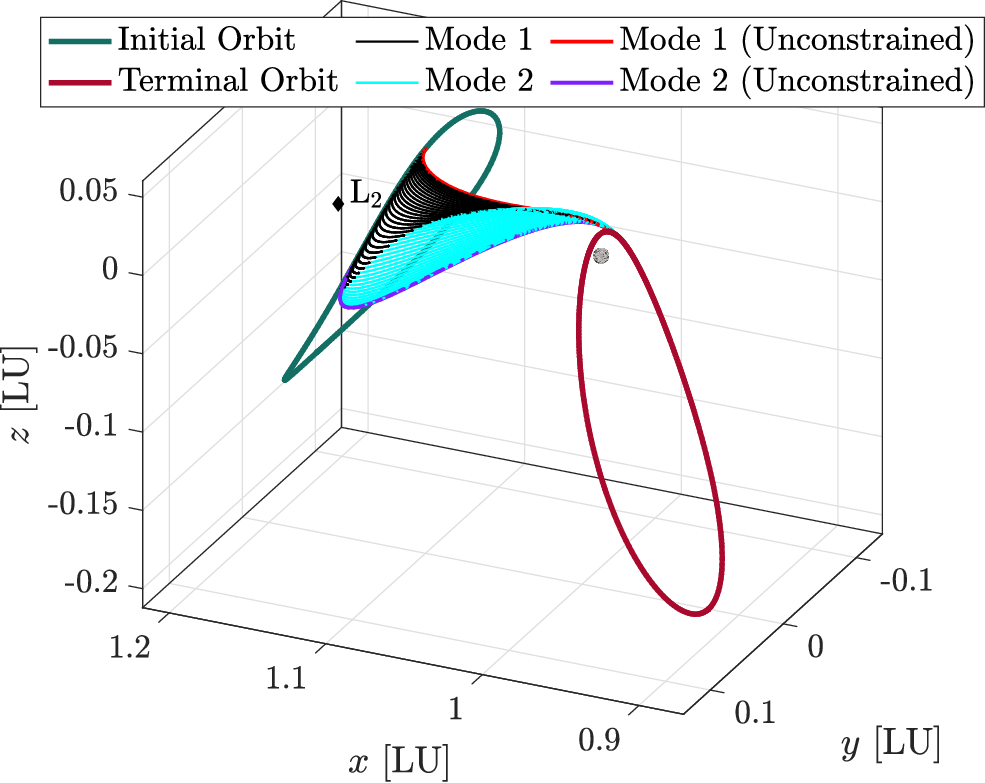}%
            \label{fig:res-mm-ts2_3d}
            }
        \\
        \hspace*{\fill}%
        \subfloat[$xy$-plane.]{%
            \includegraphics[height=5.3em]{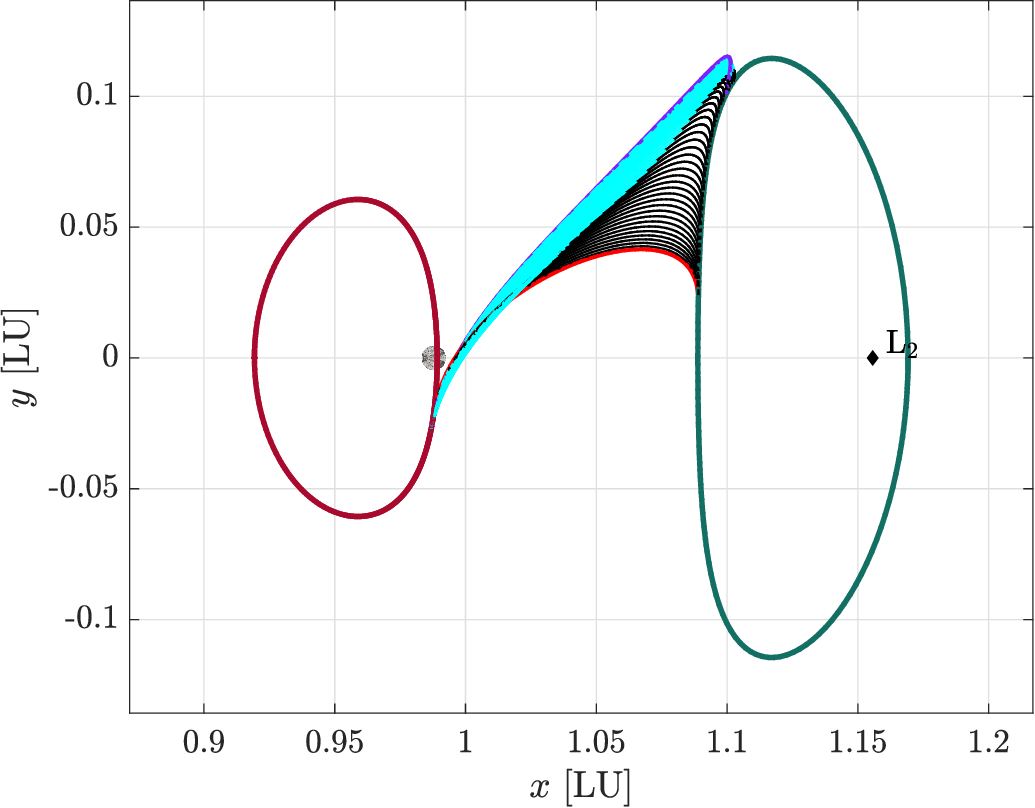}%
            \label{fig:res-mm-ts2_xy}
            }
        \subfloat[$xz$-plane.]{%
            \includegraphics[height=5.3em]{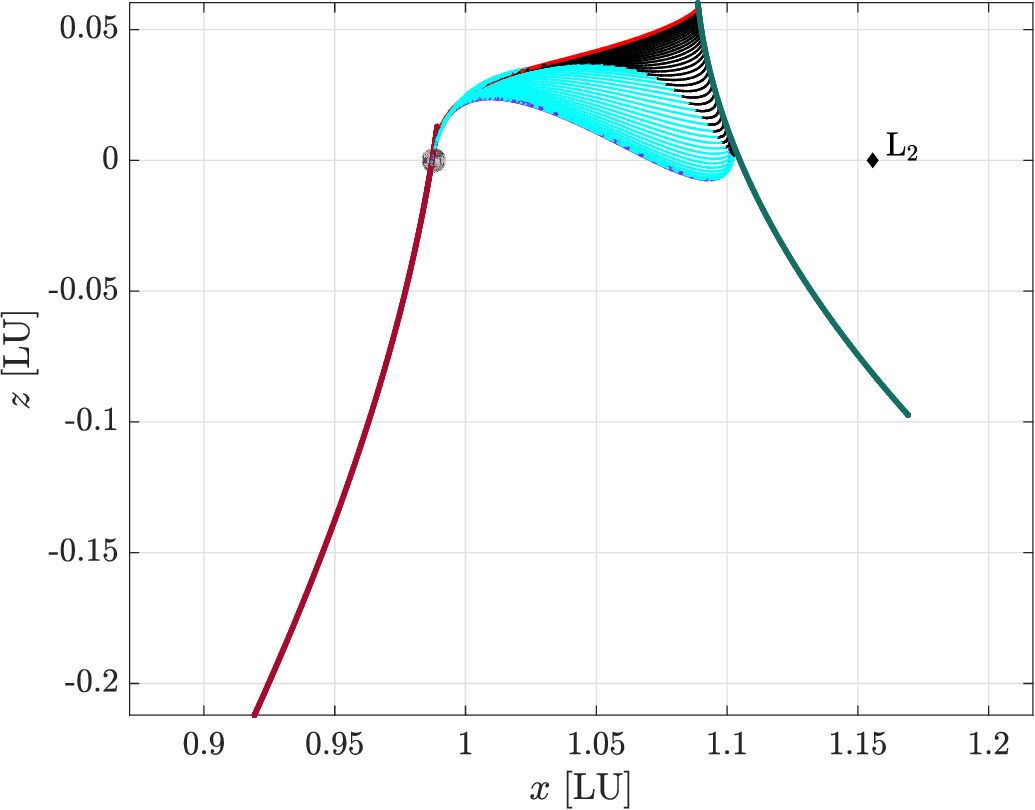}%
            \label{fig:res-mm-ts2_xz}
            }
        \subfloat[$yz$-plane.]{%
            \includegraphics[height=5.3em]{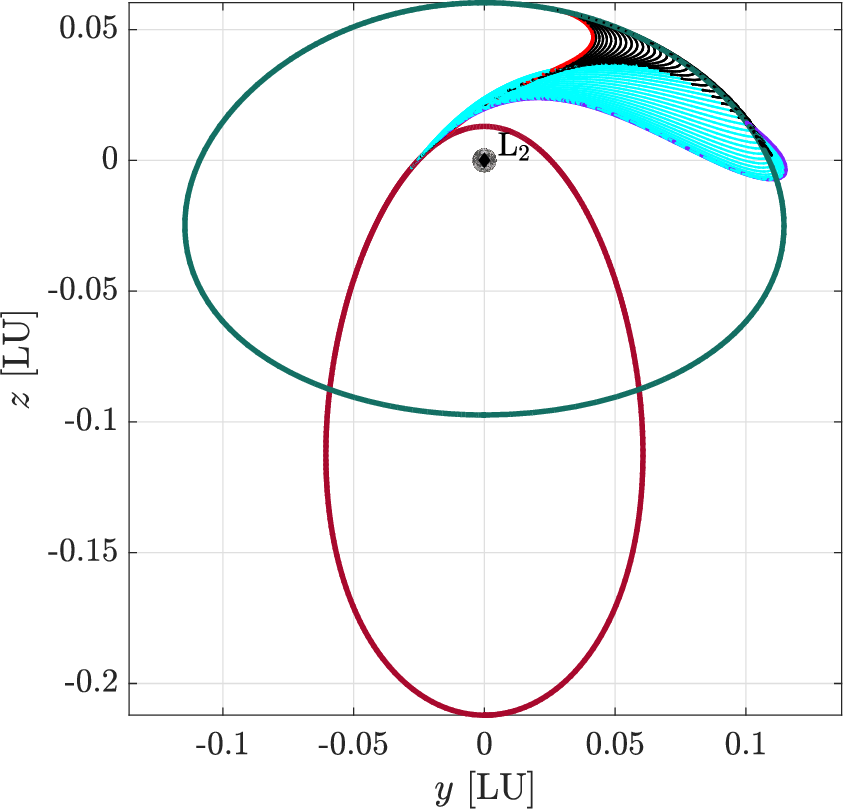}%
            \label{fig:res-mm-ts2_yz}
            }
        \hspace*{\fill}%
        \caption{Propellant-constrained, minimum-time transfers using multi-mode propulsion ($T_{\max,2}=0.25$ [N]).}
        \label{fig:res-mm-ts2}
    \end{minipage}
\end{figure}

The three-phase propellant-constrained, minimum-time problem is solved for propellant constraint values of 40 [kg] and less using multi-mode propulsion with the optimal control structure determined in Section~\ref{subsec:res-base}.
As a result, the transfer phase is partitioned into 3 domains to optimize the throttle switch times.
The throttle structure for each domain is as follows: (1) mode 1: on, mode 2: off, (2) mode 1: off, mode 2: on, and (3) mode 1: on, mode 2: off.
Note that an active mode is equivalent to full throttle.
Using the numerical approach of Section~\ref{sec:num}, the trajectory envelopes using the spacecraft specifications of Tables~\ref{tab:phys-sc-ts1} and~\ref{tab:phys-sc-ts2} are shown in Figures~\ref{fig:res-mm-ts1} and~\ref{fig:res-mm-ts2}, respectively, where the mode 1 thrusting regions are highlighted in black and the mode 2 thrusting regions are highlighted in cyan.
For reference, the propellant-unconstrained mode 1 only solution shown in Figure~\ref{fig:res-base-m1f} is highlighted in red in Figures~\ref{fig:res-mm-ts1} and~\ref{fig:res-mm-ts2}.
Relative to the trajectory envelope shown in Figure~\ref{fig:res-m1-m1} using only mode 1 propulsion, the trajectory envelopes shown in Figures~\ref{fig:res-mm-ts1} and~\ref{fig:res-mm-ts2} using multi-mode propulsion are much tighter.

For the results shown in Figures~\ref{fig:res-mm-ts1} and~\ref{fig:res-mm-ts2}, which correspond to trajectories obtained using $T_{\max,2}=0.5$ [N] and $T_{\max,2}=0.25$ [N], respectively, propellant constraint values from 40 [kg] to 1 [kg] are analyzed, where the 40 [kg] propellant-constrained trajectories are closest to the propellant-unconstrained solution using only mode 1.
For the 40 [kg] propellant-constrained solution with $T_{\max,2}=0.5$ [N], an optimal objective value of $\mcl{J}^{*}=0.287961$ (that is, 1.169 [days]) is obtained.
The optimal control profile consists of the following: (1) a 1.080 [days] mode 1 thrusting arc consuming 38.066 [kg] of fuel, (2) a 0.034 [days] mode 2 thrusting arc consuming 0.048 [kg] of fuel, and (3) a 0.055 [days] mode 1 thrusting arc consuming 1.934 [kg] of fuel.
For this solution, 40.048 [kg] of fuel is consumed in total.
It is important to note again that the propellant constraint is placed only on the fuel consumed by mode 1.
On the other hand for the 40 [kg] propellant-constrained solution with $T_{\max,2}=0.25$ [N], an optimal objective value of $\mcl{J}^{*}=0.287922$ (that is, 1.173 [days]) is obtained.
The optimal control profile consists of the following: (1) a 1.070 [days] mode 1 thrusting arc consuming 37.692 [kg] of fuel, (2) a 0.038 [days] mode 2 thrusting arc consuming 0.027 [kg] of fuel, and (3) a 0.065 [days] mode 1 thrusting arc consuming 2.308 [kg] of fuel.
For this solution, 40.027 [kg] of fuel is consumed in total.
Minimum-time transfers using only mode 2 propulsion are also obtained and highlighted in purple in Figures~\ref{fig:res-mm-ts1} and~\ref{fig:res-mm-ts2}, where optimal objective values of $\mcl{J}^{*}=0.454345$ (that is, 1.850 [days]) and $\mcl{J}^{*}=0.674895$ (that is, 2.749 [days]) are obtained when using $T_{\max,2}=0.5$ [N] and $T_{\max,2}=0.25$ [N], respectively.
For these solutions, mode 2 consumes 2.629 [kg] and 1.953 [kg] of fuel in total, respectively.
As expected, an increase in transfer time is observed when the maximum thrust of mode 2 is decreased.

\begin{figure}[!t]
    \centering
    \begin{minipage}{0.48\textwidth}
        \centering
        \hspace*{\fill}%
        \includegraphics[height=16.0em]{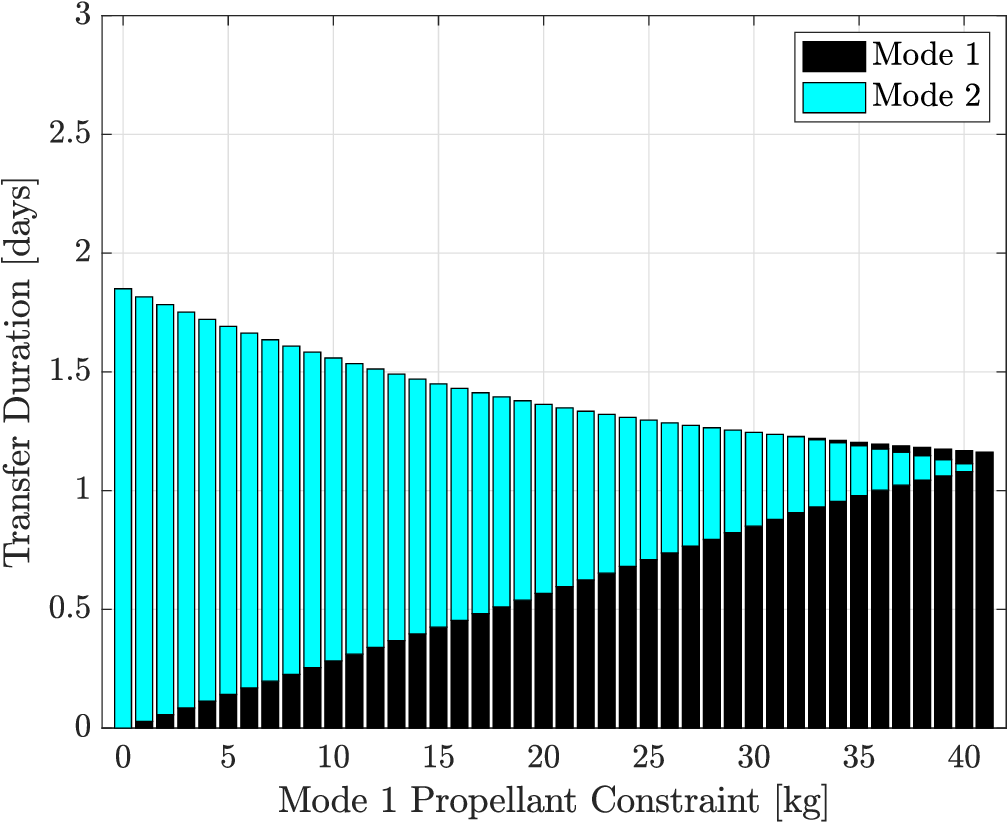}%
        \hspace*{\fill}%
        \caption{Transfer phase thrusting history for transfers shown in Figure~\ref{fig:res-mm-ts1}.}
        \label{fig:res-mm-tpdursplit1}
    \end{minipage}
    \hspace*{\fill}%
    \begin{minipage}{0.48\textwidth}
        \centering
        \hspace*{\fill}%
        \includegraphics[height=16.0em]{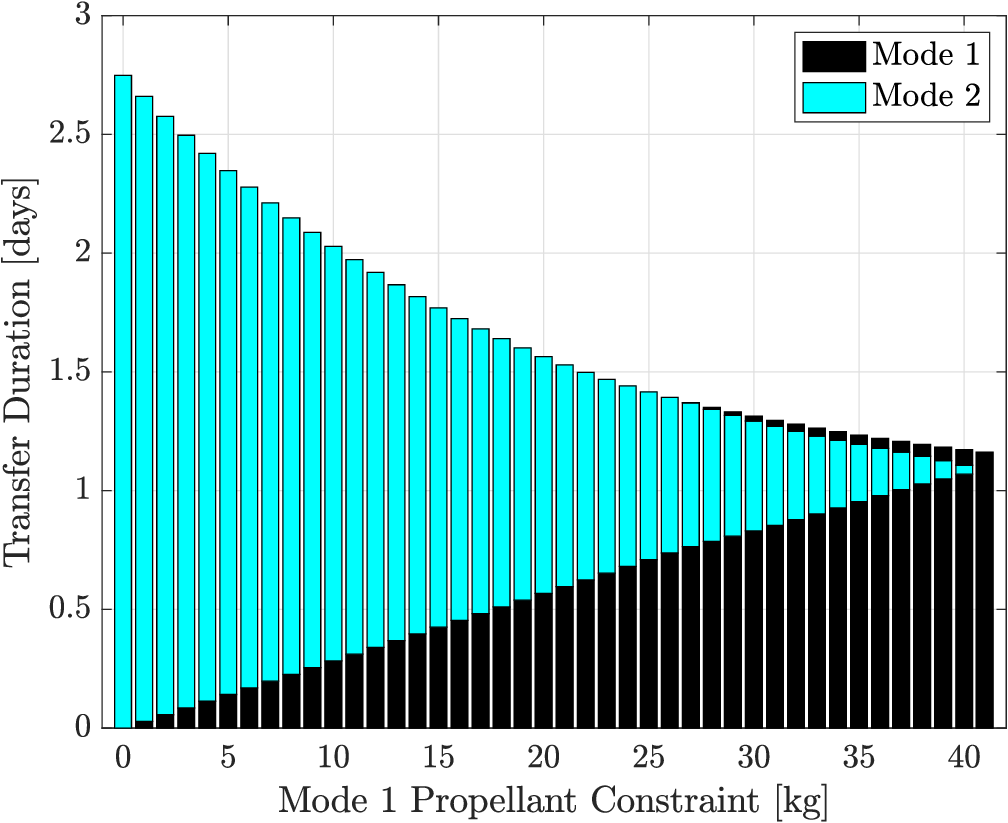}%
        \hspace*{\fill}%
        \caption{Transfer phase thrusting history for transfers shown in Figure~\ref{fig:res-mm-ts2}.}
        \label{fig:res-mm-tpdursplit2}
    \end{minipage}
\end{figure}

A summary of the optimal multi-mode control profiles for the transfer phases corresponding to Figures~\ref{fig:res-mm-ts1} and~\ref{fig:res-mm-ts2} are shown in Figures~\ref{fig:res-mm-tpdursplit1} and~\ref{fig:res-mm-tpdursplit2}, respectively.
As mentioned in Section~\ref{subsec:res-base}, the propellant-unconstrained, minimum-time transfer using only mode 1 propulsion consumes 40.973 [kg] of fuel, which is shown on the far right in Figures~\ref{fig:res-mm-tpdursplit1} and~\ref{fig:res-mm-tpdursplit2} with no mode 2 thrusting arc; furthermore, the minimum-time transfers using only mode 2 propulsion are shown on the far left with no mode 1 thrusting arc.
Relative to Figure~\ref{fig:res-m1-tpdursplit}, similar trends are observed in Figures~\ref{fig:res-mm-tpdursplit1} and~\ref{fig:res-mm-tpdursplit2}, where the duration of the mode 1 thrusting arc departing from the initial orbit decreases as the propellant constraint value is decreased.
In addition, the coast arcs are replaced with mode 2 thrusting arcs.
For $T_{\max}=0.5$ [N] with propellant constraint values ranging from 40 [kg] to 33 [kg], the optimal mode 1 and 2 throttle structures of on-off-on and off-on-off, respectively, are observed, which correspond to two switches in the control.
However, for propellant constraint values from 32 [kg] to 1 [kg], the duration of the second mode 1 thrusting arc converges to zero (that is, the corresponding thrusting arc is removed from the solution by the optimizer).
As a result, the corresponding optimal mode 1 and 2 control profiles become on-off and off-on, respectively, which correspond to only one switch in the control as shown in Figure~\ref{fig:res-mm-tpdursplit1}.
A similar phenomenon is encountered for $T_{\max}=0.25$ [N].
For propellant constraint values ranging from 40 [kg] to 27 [kg], two switches in the control are observed; however, for propellant constraint values ranging from 26 [kg] to 1 [kg], only one switch in the control is observed, which is best highlighted by Figure~\ref{fig:res-mm-tpdursplit2}.
Given these trends for the same maximum mode 1 thrust value, it is expected that as the maximum mode 2 thrust value is decreased, the propellant constraint value in which the second control switch disappears also decreases.
Finally, it is important to note that convergence issues are not encountered for either multi-mode case with any of the propellant constraint values analyzed.


\subsection{Summary and Discussion}\label{subsec:res-sum}

\begin{figure}[!b]
    \centering
    \hspace*{\fill}%
    \includegraphics[height=16.0em]{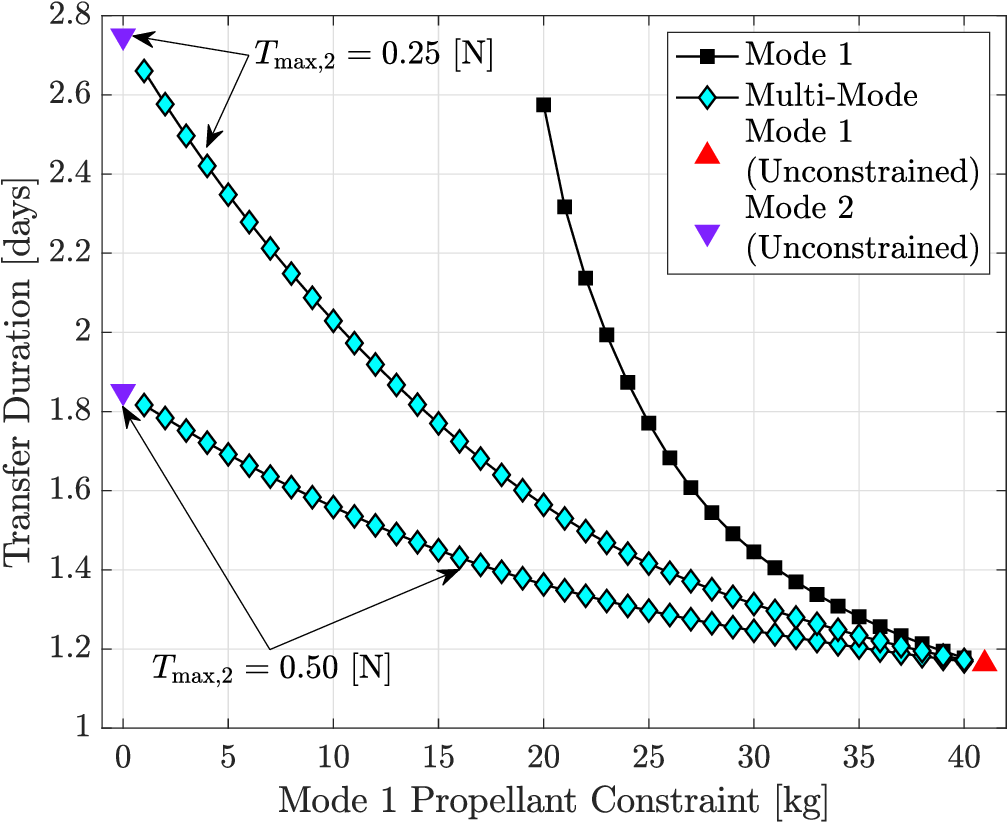}%
    \hspace*{\fill}%
    \caption{Summary of optimal transfer durations as a function of propellant constraint values.}
    \label{fig:res-sum-tdurhistory}
\end{figure}

In this section, the results are summarized for the propellant-constrained, minimum-time transfers using only mode 1 and multi-mode propulsion as shown in Sections~\ref{subsec:res-m1} and~\ref{subsec:res-mm}, respectively.
A summary of the optimal transfer duration as a function of the propellant constraint value is shown in Figure~\ref{fig:res-sum-tdurhistory}.
For the propellant-unconstrained, minimum-time transfer using only mode 1 propulsion, an optimal transfer duration of 1.163 [days] is obtained, where 40.973 [kg] of fuel is consumed.
For propellant constraint values of 40 [kg] to 20 [kg] using only mode 1 propulsion, the optimal transfer durations range from 1.178 [days] to 2.575 [days], and the corresponding total fuel consumption ranges from 40 [kg] to 20 [kg], respectively.
For propellant constraint values of 40 [kg] to 1 [kg] using multi-mode propulsion with $T_{\max,2}=0.5$ [N], the optimal transfer durations range from 1.169 [days] to 1.816 [days], and the corresponding total fuel consumption ranges from 40.048 [kg] to 3.541 [kg], respectively.
For the corresponding propellant-unconstrained, minimum-time transfer using only mode 2 propulsion, an optimal transfer duration of 1.850 [days] is obtained, where 2.629 [kg] of fuel is consumed.
For propellant constraint values of 40 [kg] to 1 [kg] using multi-mode propulsion with $T_{\max,2}=0.25$ [N], the optimal transfer durations range from 1.173 [days] to 2.660 [days], and the corresponding total fuel consumption ranges from 40.027 [kg] to 2.870 [kg], respectively.
For the corresponding propellant-unconstrained, minimum-time transfer using only mode 2 propulsion, an optimal transfer duration of 2.749 [days] is obtained, where 1.953 [kg] of fuel is consumed.
For all results previously discussed, mode 1 consumes the maximum allowable amount of fuel when utilizing either only mode 1 propulsion or multi-mode propulsion.

Four key features in the evolution of the transfers are highlighted in Figure~\ref{fig:res-sum-tdurhistory}, which demonstrate the trade-off between transfer duration and propellant consumption for single- and multi-mode propulsion systems.
First, for a given propellant constraint value, using multi-mode propulsion decreases the transfer duration relative to its single-mode counterpart, with the reduction in transfer duration significantly increasing as the propellant constraint value is decreased.
Second, for a given transfer duration, using multi-mode propulsion decreases the required propellant amount relative to its single-mode counterpart.
Third, as the magnitude of the maximum thrust of mode 2 is decreased, the transfer time increases, as expected; however, the control structure still changes depending on the propellant constraint value, where a different number of control switches is observed for lower propellant constraint values when compared to the single-mode solution.
Finally, using multi-mode propulsion increases the overall feasibility of the transfer given various propellant constraint values, which is otherwise limited by the minimum-fuel solution obtained using single-mode propulsion.



\section{Conclusions}\label{sec:conc}

A numerical trajectory optimization study is conducted for propellant-constrained, minimum-time transfers between periodic orbits in the Earth–Moon elliptic restricted three-body problem using both single- and multi-mode propulsion systems.
A novel three-phase optimal control problem is formulated in which the true anomaly on the primary orbit is employed as the independent variable, which eliminates the need to solve Kepler’s equation at arbitrary epochs.
To facilitate evaluation of the minimum-time objective functional, an additional state variable is introduced through normalized time, while a path constraint ensures that only one propulsion mode may be active at any given epoch.
The three phases consist of the following: (1) a coast phase along the initial periodic orbit, (2) a controlled transfer phase, and (3) a coast phase along the terminal periodic orbit.
Then, the controlled transfer phase is further partitioned into domains in order to optimize the switches in the throttles of each mode.
The resulting optimal control problem is solved numerically using an adaptive Gaussian quadrature direct collocation method.
Case studies are performed for transfers from an $\tL_2$ southern halo orbit to a near-rectilinear halo orbit across multiple single- and multi-mode propulsion architectures and a wide range of propellant constraint values.
Results demonstrate that multi-mode propulsion reduces the transfer time for a fixed propellant budget and lowers propellant requirements for a fixed transfer time relative to single-mode systems.
Moreover, multi-mode throttle structures often exhibit a different number of switches compared to their single-mode counterparts, and the use of multi-mode propulsion increases transfer feasibility in regimes where single-mode solutions are otherwise infeasible.


\section*{Acknowledgements}

The authors gratefully acknowledge support for this research from the U.S. National Science Foundation under Grant CMMI-2031213 and the NASA Florida Space Grant Consortium under Grant 80NSSC20M0093.


\renewcommand{\baselinestretch}{1}
\normalsize\normalfont



\end{document}